\documentclass[11pt]{article}

\usepackage{lineno}
\usepackage{amsmath}
\usepackage{amsfonts}
\usepackage{adjustbox}
\usepackage{siunitx}
\usepackage{booktabs}
\usepackage{float}
\usepackage{multirow}
\usepackage{color}
\usepackage{textcomp}
\usepackage{comment}
\usepackage{leftidx}
\usepackage{natbib}
\usepackage{pifont}% http://ctan.org/pkg/pifont
\newcommand{\cmark}{\ding{51}}%
\newcommand{\xmark}{\ding{55}}%

\usepackage{hyperref}
\usepackage{array}
\newcolumntype{-}{>{\global\let\currentrowstyle\relax}}
\newcolumntype{^}{>{\currentrowstyle}}

\newcommand{\scc}[0]{{\textrm{sc}}}
\newcommand{\fc}[0]{{\textrm{fc}}}
\newcommand{\cc}[0]{{\textrm{cc}}}

\newcommand{\red}[1]{{\color{black}#1}}
\newcommand{\blue}[1]{{\color{black}#1}}
\newcommand{\ObjVal}{\mathrm{ObjVal}} 

\title{Three Network Design Problems for Community Energy Storage}

 \author{
Bissan Ghaddar\thanks{Ivey Business School, Western University, \href{mailto:bghaddar@ivey.ca}{bghaddar@ivey.ca}}
\and Ivana Ljubi\'c\thanks{ESSEC Business School of Paris, France, \href{mailto:ivana.ljubic@essec.edu}{ivana.ljubic@essec.edu}} \and Yuying Qiu\thanks{ESSEC Business School of Paris, France,  \href{mailto:yuying.qiu@essec.edu}{yuying.qiu@essec.edu}}
}

\date{}
\begin{document}
\maketitle
\begin{abstract}
In this paper, we develop novel mathematical models 
to optimize utilization of
%distributed energy resources and 
\red{community energy storage (CES)} by clustering prosumers and consumers into energy sharing communities/microgrids in the context of a smart city.
Three different microgrid configurations are modeled 
using a unifying mixed-integer linear programming formulation. These configurations \red{represent three different business models, namely: }
%modeling the microgrids using: 
the island model, the interconnected model, and the Energy Service Companies model. 
%for optimal clustering of communities and the corresponding CESs' configuration in terms of size and location is so far lacking. 
\red{The proposed mathematical formulations determine the optimal households' aggregation as well as the 
location and sizing of CES. }
To overcome the computational challenges of treating operational decisions within a multi-period decision making framework, we \red{also} propose a decomposition approach to accelerate the computational time needed to solve larger instances. 
%We provide the optimal community energy storage locations and capacities 
\red{We conduct a case study based on real power consumption, power generation, and location network data from Cambridge, MA}. Our \red{mathematical models} and the underlying algorithmic framework can be used in operational and strategic planning studies on smart grids to incentivize the communitarian distributed renewable energy generation and to improve the self-consumption and self-sufficiency of the energy sharing community. The models are also targeted to policymakers of smart cities, utility companies, and Energy Service Companies as the proposed models support decision making on renewable energy related projects investments.
\end{abstract}
\noindent
Keywords: OR in energy; Energy management system; Community energy storage; %Self-consumption \sep
Decentralized Energy Resources; Mixed integer linear programming %\sep Benders decomposition

\section{Introduction}
As we move into a sharing society and smart cities' structure, energy sharing within a neighborhood will become more common thanks to the development of new technologies for energy generation and storage as well as smart grids. Households can trade energy produced from solar panels with each other on the local energy market via community microgrids, or sell the surplus to the utility companies via the wider power grid. \red{Participants within a local energy community often share the costs and benefits of renewable energy projects, storage systems, or other energy-related initiatives. They may also collectively manage and trade surplus energy among themselves contributing to an increased energy sustainability. They can be both grid-connected and off-grid, and they may integrate with existing utility infrastructure or operate independently, depending on local regulations and the community's goals. } %Smart grids therefore play an important role in optimizing resource allocation efficiently as well as shifting the traditional fossil fuels consumption to clean renewable energy. 

Many countries have experienced an increase in ``renewable energy communities'', in which groups of households and other entities within a neighborhood are motivated to reduce their energy costs and promote the deployment of renewable energy \citep{doci2015let}. Concurrent with the increasing renewable energy penetration in distributed smart grids, more and more households are inclined to install roof-top photo-voltaic (PV) systems to lower their electricity bills. Therefore, the traditional energy systems and energy grids are being challenged, as the new decentralized energy community %will 
tend to become one of the main \red{promoters of} additional renewable energy penetration. To better cope with the imbalance between supply and demand caused by the increasing intermittent renewable energy penetration in smart grids, a new energy network design focusing on decentralized energy systems is an effective approach in smart grids {\citep{alstone2015}}.

In this article, we study the design of \red{local energy communities} using community energy storage (CES)%\footnote{A summary of abbreviations used in this paper is presented in the Appendix}, 
as a possible alternative to single household batteries. In a community-based prosumer setup, the households equipped with solar panels act as their own energy providers and can exchange electricity surplus with other households within the community. This network connects prosumers and traditional consumers in a given neighborhood, so that the produced energy is dispatched through a central energy management system (EMS) in a coordinated way. Community energy storage is a crucial component of the EMS. The design of smart grids in the future will take advantage of CES in dealing with more dynamic loads and energy sources \citep{roberts2011role}. The energy stored in the CESs guarantees the electricity supply during peak hours, which reduces the electricity bills from utility companies of the consumers, thus bringing economic benefits to the households as well as to the community.% who may act as a single energy consumption agent. 

Energy storage can help integrate local renewable generation into existing power systems, but the questions on how to deploy the batteries within a community network to maximize the profit of the CES investment, and how to optimally dispatch the energy in the system to minimize the electricity bill of the community remain open. Selecting optimal locations and capacities for batteries at certain nodes of a given neighborhood is a complex combinatorial optimization problem. In the context of smart cities, the problem may involve simultaneous decisions for hundreds of microgrids. According to \cite{parag2016electricity}, there are three types of prosumer markets: 1) island microgrids, 2) interconnected microgrids, and 3) energy service company model. In the first one, the microgrid is isolated and operates in the so-called island model, whereas in the second and third one, the microgrid operates in connection to the main grid. In this article, we address all three types of prosumer markets, while changing the objective functions according to the incentives of the underlying \red{business model.} More precisely, in the island model, the goal is to find the optimal investment in battery storage that guarantees a seamless energy supply to all households (and any excess generation cannot be monetized unless it is stored for later usage). 
\red{In the second business model,} the purpose of each community is to reduce the electricity bill as a whole by increasing the electricity exportation revenue and decreasing the electricity consumption cost from the main grid, subject to certain battery investment constraints. In the third %setting, 
\red{business model,} there is an Energy Service Company (ESCO) that makes an upfront investment into the infrastructure and then collects a service fee that is proportional to the energy that each household gets from discharging CESs. In this case, the objective of an ESCO is to maximize the discounting cash flows over various years. There also exists the fourth possibility to have a peer-to-peer exchange network, like for example the Sonnen Community \citep{Sonnen}. This type of energy exchange poses additional operational challenges to the grid, because of the capacity bottlenecks and the amount of energy exchanged. In this paper, we focus on geographically clustered microgrids, where the renewable capacity and energy storage can be %economically 
locally clustered to improve the performance of the energy sharing community, thus peer-to-peer models are not considered in our network design. 

\subsection{Related work {and contribution}}
Energy allocation optimization with CES has gained increased attention in recent years. In \cite{alskaif2017reputation}, a CES with a centralized EMS is proposed. The EMS manages the CES and controls the allocation of available energy in the storage unit to households. The optimal operation of the community energy storage system for PV energy time-shift \citep{parra2015optimum}, demand load shifting \citep{parra2016optimum,terlouw2019multi} and some other benefits such as economies of scale, energy trading and enhanced grid balancing capabilities \citep{arghandeh2014economic} are demonstrated. Some stochastic features of the CES operations are also considered \red{in the literature.} A stochastic smart charging framework for CES in residential microgrids is proposed in \cite{alskaif2017smart}, and a day-ahead scheduling model is built in \cite{liu2018energy} to increase the operation profit of the energy-sharing network, considering the uncertainty of PV energy, electricity prices, and prosumers' load. Under the context of uncertainty of energy, a robust optimization model for long-term energy planning is proposed in \cite{moret}, which features uncertain inputs in the energy planning practice.

The optimal energy management for energy storage systems has been frequently studied.  %recently. %, but few publications are found in the classical management and decision science literature transferring knowledge from well-known fields of research, such as inventory control, to the discussed application of energy management. 
A systematic review of the literature on energy management for energy storage applications including optimization methods is conducted in \cite{weitzel}. \cite{kuznia} developed a Benders’ decomposition algorithm for a comprehensive hybrid power system design problem, including renewable energy generation, storage devices, transmission networks, and thermal generators, for remote areas using a mixed integer programming model (MILP). To determine the optimal battery configuration, \cite{chen2011sizing} conducts a cost-benefit analysis for the optimal size of an energy storage system for both the grid-connected and island model network, whereas the optimal sites and size of energy storage systems to perform spatio-temporal energy arbitrage most effectively has been identified in \cite{fernandez2016optimal}. \cite{qi2015} proposed models of transmission network planning with co-location of energy storage systems (ESS). Their models determine the sizes and sites of ESS as well as the associated topology and capacity of the transmission network under the feed-in-tariff policy instrument. This is similar to our work but on a transmission level. In \cite{sardi2017strategic}, locations, sizes and operational characteristics of CES are optimized to enhance network performances including the CES integration. \cite{Dai2018} provide optimal operational strategies for using shared storage in buildings. \cite{van2018} modeled the energy management problem between the generation companies owning centralized assets and the microgrids using a bi-level stochastic mixed integer program.

Clustering networks within smart grids is also discussed in the literature. The optimal planning of the interconnected network of multimicrogrids is presented in \cite{che2015optimal} and the microgrid clustering problem is solved by studying the optimal power flow between clusters while managing congestion and power losses in \cite{boroojeni2016novel}. Long-term planning of electric power infrastructure to choose the optimal investment strategy and operation schedule is described in \cite{lara2018deterministic}, but the authors did not consider the %siting
problem of locating storage in the residential area. In \cite{van2018techno,barbour2018community}, a systematic comparison of batteries for individual dwellings and communities is conducted. The authors provide insights into the optimal aggregation level of storage deployment and when combined with Demand Side Management (DSM), it can improve self-consumption \cite{van2018techno}. In \cite{barbour2018community}, the authors use a data-driven approach to group households into local energy sharing communities with a single CES, and they illustrate the advantages of CES compared to the household energy storage (HES), including economies of scale for batteries and benefits related to the lower likelihood of consumption peaks. 
\red{Table \ref{tab:lit} summarizes the related literature on optimization problems related to energy storage with models looking into clustering, battery sizing, battery locations, and battery scheduling. }

\begin{table}[h]
\centering \red{
\begin{tabular}{p{0.1\linewidth} | p{0.11\linewidth} | p{0.14\linewidth}  p{0.06\linewidth}  p{0.06\linewidth}  p{0.06\linewidth} | p{0.2\linewidth} }
\toprule
\multicolumn{1}{l|}{ 
}                               & Decision & \multicolumn{4}{|c|}{ 
 Modeling} &Methodology\\
 &Type &    Clustering & Bat. size & Bat. loc. & Bat. sch.& used\\
\midrule
\multicolumn{1}{l|}{\cite{che2015optimal}}                                              &strategic           & \cmark (microgrids)      & \xmark   &  \xmark & \xmark & probabilistic minimal cut set \\
 \multicolumn{1}{l|}{\cite{boroojeni2016novel}}                  &operational           & \cmark (microgrids)      & \xmark   &  \xmark & \cmark & network optimization \& power routing \\
 \multicolumn{1}{l|}{\cite{lara2018deterministic}}                                    &strategic           & \xmark      & \xmark   &  \xmark & \xmark &MILP \\
 \multicolumn{1}{l|}{\cite{van2018techno}}                   &operational           & \xmark      & \xmark   &  \xmark & \cmark &MILP \\
 \multicolumn{1}{l|}{\cite{barbour2018community}}                  &  operational         & \xmark      & \cmark   &  \xmark & \xmark & heuristics\\
 \multicolumn{1}{l|}{\cite{chang2022}}                  &     operational      & \cmark      & \xmark   &  \xmark & \cmark & MILP\\
 \bottomrule
\multicolumn{1}{l|}{This work}  &strategic           & \cmark      & \cmark   &  \cmark & \xmark &MILP with Benders \\                                     
\bottomrule
\end{tabular}}
\caption{Existing literature on shared energy storage.}
\label{tab:lit}
\end{table}

{Our models have their roots in the (capacitated) facility location problems, which have been extensively studied in the operations research literature in the last decades (see, e.g.,  books by \cite{laporte2019book,drezner2004facility,nickel2006location}). In particular, the concept of finding the optimal locations and configurations of batteries for the CES and the ``clustering'' of households, can be seen as a capacitated facility location problem with modular facilities (see, e.g., \cite{ALARCONGERBIER2022108734}) where CES locations correspond to open facilities, and the assignment of households to CES represents the desired clustering. We point out that our models are much more complex compared to the related facility location and network design problems, since the load data, power flows and state-of-charge of batteries over the planning horizon need to be taken into account. }

%Although 
\paragraph{{Contribution:}}

\red{The purpose of this paper is twofold, the first is to formulate and present different business models for community energy storage. The second is to apply these models in practice for a realistic use case and use the solutions to analyze different components of the system. This work provides a first step towards having practical models for shared energy storage.}

\begin{itemize}
    \item 
Related to the first part, even though there has been extensive research on energy storage systems in the recent literature, %however 
a mathematical model %of the
for optimal clustering of communities and the corresponding CESs' configuration in terms of size and location is so far lacking. Our research aims at closing this gap by providing an exact algorithm to get the optimal prosumers' aggregation as well as the %siting 
location and sizing of community %level of 
energy storage. {To this end, we propose a generic modeling framework and new mathematical models to deal with three different types of microgrids, \red{representing three different business models based on } 
%namely, 
the island microgrids, interconnected microgrids, and the third-party energy service company model.  Clustering thousands of households poses a serious challenge for the modern state-of-the-art mixed integer programming solvers. Benders decomposition approaches have been shown to provide some of the state-of-the-art results for the related facility location and network design problems (see, e.g., \cite{FischettiLS16,FischettiLS17,CordeauFL19,ConiglioFL22,Duran-Mateluna-et-al:2022,Gaar-Sinnl:2022}). This has also motivated our algorithmic choice: rather than giving the compact models to an off-the-shelf solver, we aim at assessing the ability of Benders decomposition in improving the computational performance and finding optimal or high-quality solutions for larger instances of realistic size. 

In this article, we assess the tractability of proposed MILP models when it comes to solving instances with realistic input data. 
\red{We aim to answer the following questions: 
\begin{itemize}
    \item how big are the instances that can be solved to optimality (by using the compact models, or Benders decomposition instead)?
    \item how large the gaps are for the instances where the proof of optimality cannot be achieved? 
\end{itemize}
\red{In this article we study deterministic mathematical models, starting from given consumers' and prosumers' load generation profiles. We develop models using the deterministic rather than the uncertain data, due to the problem complexity inherited by their facility location counterparts.  All three mathematical models are new, and as it is common in the Operations Research literature, before studying less tractable stochastic optimization models, we aim to understand the structure and complexity of their deterministic counterparts. }
These first experiments help us to understand where are the limits of the proposed MILP models. 
They also serve as guidance for finding where to put further modeling efforts in order to deal with data uncertainty. }
%To assess the effects of data uncertainty on our models, we provide a-posteriori analysis where we vary the demand data for different days to see the impact on the operational costs for the three proposed models.  % \red{TODO: complete}

\item \red{For the second contribution, we use a realistic data set derived from Cambridge, MA dataset of \citet{barbour2018community},  where real power consumption and generation data of households is considered within a given community. We apply the three different business models to the resulting dataset. %Sensitivity analysis is performed to show the impact of the variation in the demand and production data on the operational costs of the community. 
To assess the effects of data uncertainty on our models, we provide a posteriori analysis where we vary the demand data to see the impact on the operational costs of the community for the three proposed models. Different business models might be more attractive for different types of communities based on consumption, preferences, and the number of prosumers.}
}
\end{itemize}

The rest of the paper is organized as follows. The three types of 
%clustering 
\red{business models} as well as %how we
the creation of the local communities are presented in Section \ref{sec:3models}, the corresponding microgrid models as well as the mathematical models and the Benders decomposition approach are described in Section \ref{sec:formulation}. A use case using real data from Cambridge, MA is presented and analysed in Section \ref{sec:usecase}. Numerical results are performed in Section \ref{sec:results}. Finally, Section \ref{sec:conclusion} provides the main conclusions and discusses potential 
future research directions.

\section{Three Microgrid Models}\label{sec:3models}
The residential energy sharing network consists of a set of households, including the ones who install the solar panels on their rooftops, which are referred to as prosumers, as well as the consumers who only get \red{electricity} from a shared battery or from utility companies. These households are connected within a community microgrid which has a star \red{topology} network where each star is centered around a single community energy storage. %\footnote{It is misleading to call it microgrid: she wanted to introduce the concept of stars, where several households are connected to a CES. We need to call it clusters, or somehow differently} 
In this work we consider that every household is equipped with a smart meter to record their electricity transaction. We therefore discuss three types of %clustering methods 
\red{business models} in the sections that follow.
\subsection{Island model}	\label{sec:IL}
\begin{figure}[H]
     \centering
\includegraphics[width=0.3\textwidth]{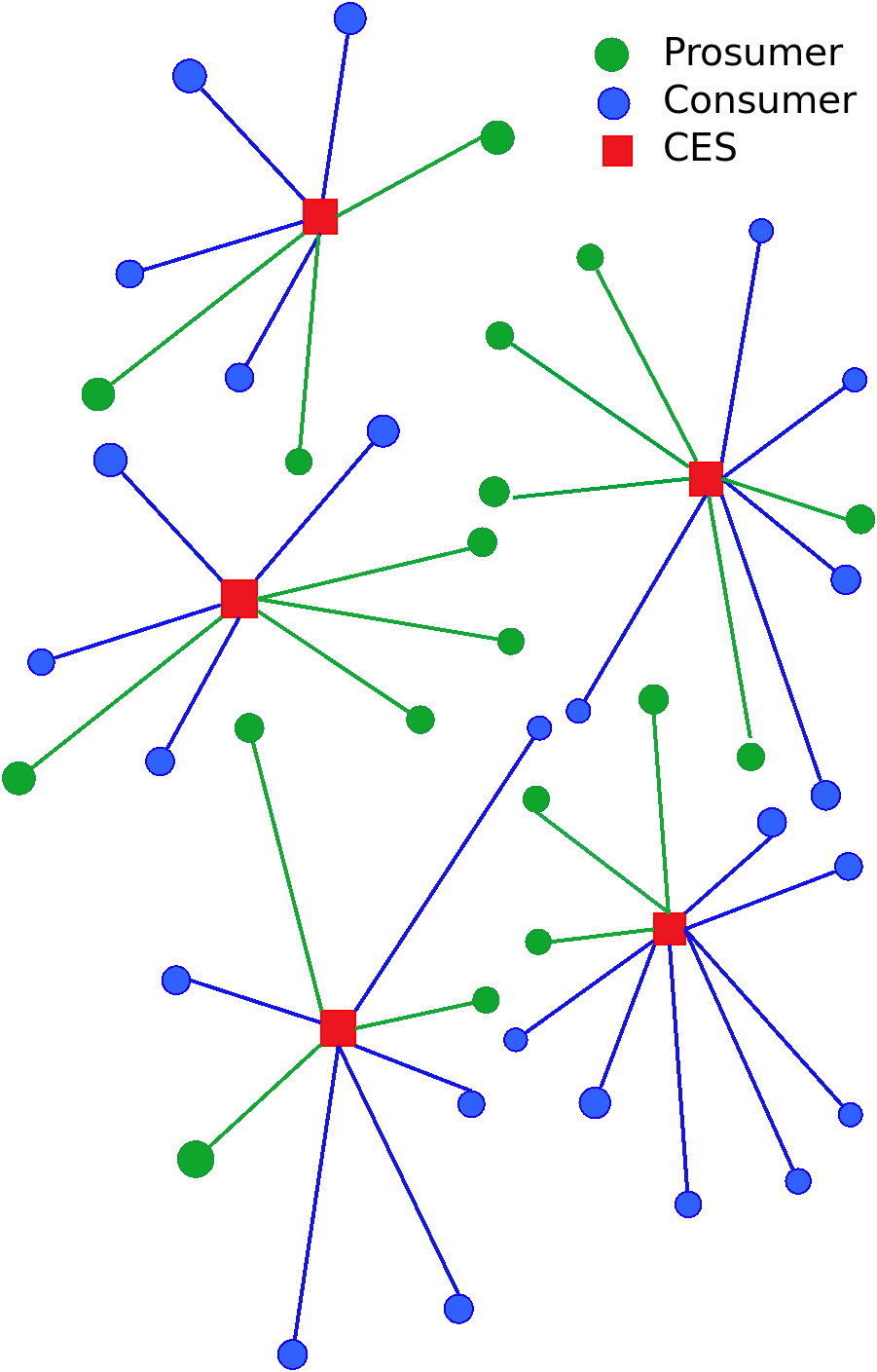}
\caption{A simple architecture of island model clusters}
\label{fig:CFL1}
\end{figure}

As a basic %model 
\red{way of clustering,} the island model allows the community microgrid to operate autonomously without any \red{electricity} purchase or export to the main grid. Therefore, it is crucial to dispatch the \red{electricity} that prosumers generate at the microgrid level via the central EMS using CESs, in order to guarantee the self-sufficiency of each household at peak times. Most island model communities exist in rural areas where the connection with the main grid is difficult. For example, ABB built a 100\% sustainable resort in the Maldives \citep{ABB}; Schneider Electric worked on transforming existing interconnected microgrid into island model micogrid such as ``Boston One Campus Islandable Microgrid'' \citep{island}; Indigo Power is moving forward with the Upper Murray Islandable Micro-Grid Project \citep{indigo}. In the island clustering model, we adopt the single bus microgrid which assumes that all prosumers, consumers, and storage devices interact at one node and thus do not generate distribution lines losses \citep{weitzel} in a radial network. \red{We assume that potential battery locations are a subset of prosumer locations as not all prosumers will have one and they are selected such that they cover the geographical area uniformly for a given community.} Figure \ref{fig:CFL1} illustrates the architecture of the island clustering model of one small neighborhood.
The CES can help optimize the prosuming services with excess generation only to the point where load shifting services are available \citep{parag2016electricity}, and the microgrid is required to supply critical loads of multiple parties to minimize the social welfare loss of the community \citep{li2016optimal}. In the island model, %since there is no external energy supply, the existence of shiftable loads and ensuring the electricity supply for critical loads is important. 
{if there is no external energy supply, 
%the existence of shiftable loads and 
ensuring the energy supply for critical loads is important. 
This can done through a 
%DER can be either 
distributed generation. In our case,  }
to avoid electricity shortage in case of insufficient generation (even with the smart load shifting), we assume that each microgrid is equipped with at least one microturbine to ensure critical load supply.   Microturbines (MT) are small single-staged combustion turbines that generate from a few kWs to a few MWs of electricity, powered by natural gas, biogas, hydrogen, or diesel. {They typically complement renewable resources to overcome the variable nature of these sources and have a self sufficient microgrid \citep{pourmousavi2010}.} %Considering the installation cost and the fuel cost, the total cost function of microturbines is $F(P_{i,mt}^t)=a\sum_{i \in B}y_{i,mt}+b \sum_{t\in T}\sum_{i\in B}P_{i,mt}^t$ where $a$ and $b$ are the installation and fuel cost coefficients respectively, $y_{i,mt}$ indicates whether a mico-turbine is installed at location $i$ or not, and $P_{i_mt}^t$ is the amount of energy produced by the microturbine located at node $i$ at period $t$. 
The transformation from the traditional energy supply network to the island microgrid based decentralized community is supposed to be undertaken by the utility companies. {Because in this scenario the utility is in charge of the microgrid network transformation and the deployment of the battery, the objective of clustering in the island model is to minimize the battery installation cost and microturbine installation and operations fuel cost to satisfy the demand from households of each community.}%\footnote{Here we claim that we do not pay installation costs for the microtourbines, but they are in the objective function of the island model. What do we have in the python code? \blue{we have the investment costs for both CES and MT in the python code, it is the same as equation (1), so I rephrased this sentence}}

{Another potential application of the island model is for \emph{intentional islanding}, in which case the microgrid is usually connected to the utility system, hence the electricity supply for critical loads is ensured \citep{Balaguer-et-al:2011}.  However, in some extreme situations when the microgrid is cut off from the main grid (due to a catastrophic utility failure or some other unpredictable external events), the microgrid is supposed to continue to provide adequate power to the load.
%describes the condition in which a microgrid or a portion of the power grid, which consists of a load and a distributed generation (DG) system, is isolated from the remainder of the utility system. In this situation, 
Hence, the clustering of the microgrid has to be made in such a way that it remains self-sufficient in case it becomes isolated from the remainder of the utility system.}

\subsection{Interconnected model}
In the interconnected model, the electricity demand for consumers is supplied by either a community battery from microgrid, or the utility company from the main grid, while the prosumers' electricity demand can be satisfied by community battery, the main grid connection, or their own rooftop PV energy generation. All the households are connected to the main grid so that the utility company can provide electricity when the prosumers of one community cannot supply enough \red{electricity} for all households of their network or there is not enough \red{electricity} stored in the CES. On the other hand, if the community generates excessive \red{electricity} and the CESs are charged to their maximum capacity, the prosumers can also export electricity to the main grid to get extra compensation from the utility company. For example ``Repowering London'' \citep{repower} provides service for community owned renewable energy project. As of March 2022, Repowering London has installed 670kWp of community owned renewable energy. {Additionally Manitoba Hydro \citep{manitoba} is using energy storage and grid-connected solar PV systems to reduce the amount of electricity bills of households}. This clustering setting can promote services of companies like ``Repowering London'' and ``Manitoba Hydro'' but at the same time, it is the most challenging for utility companies, because each microgrid works as either load or generator, thus it requires the wide power network system to be flexible and resilient to deal with bidirectional power flow. 
\begin{figure}[H]
     \centering
\includegraphics[width=0.5\textwidth]{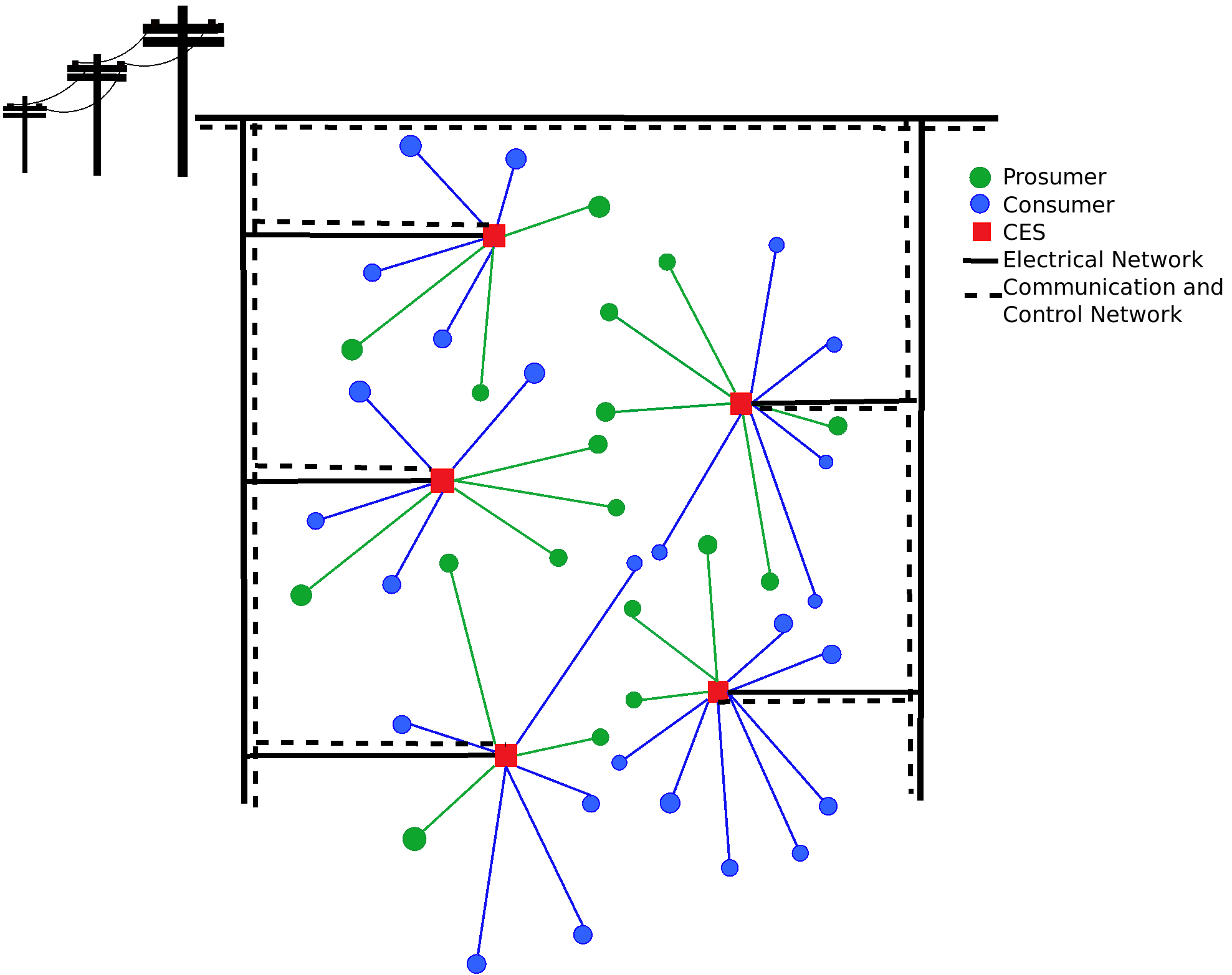}
\caption{A simple architecture of interconnected model clusters. \red{The architecture is taken from the battery storage perspective, not showing the consumers' or prosumers' connection to the grid.}}
\label{fig:CFL2}
\end{figure}
Figure \ref{fig:CFL2} illustrates the architecture of the interconnected %clustering 
model of one small neighborhood. In order to maximize the social welfare from the community clusters, the objective of clustering in this interconnected model is to maximize the Net Present Value (NPV) of net energy generation income for all the households in a given neighborhood. {To emphasize the economical effect that CES can bring to the community and to incentivize the household to install batteries for the community network, the battery deployment budget should be less than the total electricity expenses (when the community imports energy exclusively from the main grid without CES), for each community within the considered planning horizon. In this %cluster 
setting, since the prosumers can get the most benefit from the deployment of CES by storing their electricity in the battery and exporting surplus to the main grid, we assume that prosumers of each cluster are charged the CES deployment expenses.}
%\footnote{Is this consistent with our math model and with the implementation?\blue{Yes: the objective is the same as equation (6) so its only the operational costs but there is a constraint (7) that forces the investment costs to be lower than the costs just considering the grid.}}

\subsection{ESCO model}
An Energy Service Company is an energy-related commercial entity that provides energy solutions for their clients. In this setting, we assume that the ESCO takes charge of the design and implementation of energy savings projects including the CESs deployment. The ESCO starts by performing an analysis of the property, designs an energy efficient solution, installs the required elements, and maintains the system to ensure energy savings during the payback period \citep{ESCO}. Existing ESCOs such as ``Schneider Electric'', ``Siemens'', ``Ameresco'', and ``Enerpower'' can provide energy storage auxiliary services for communities along with PV solar panel installation services. As the intermediary between households and the utility company, ESCOs have two main functions in an energy sharing community. First, they guarantee to consumers and prosumers the reduction of their electricity bill by optimally allocating the energy that prosumers produce; Second, they can help the utility company to manage the unpredictable total community load or generation by increasing the self-consumption and self-sufficiency of each microgrid. The global ESCO performance contracting market is expected to grow from \$25.2 billion in 2020 to \$49.6 billion in 2030 \citep{esco_growth}.

\red{The ESCO will install and setup the battery storage and in return charge the households based on the energy they consume as well as pay the prosumers based on the amount of energy produced.} The business model for the service that ESCOs provide in each energy sharing community is that: they charge the consumers and prosumers for the energy that they use from the CESs at a relatively low fixed internal electricity buying price and pay the prosumers for the energy that they produce at a fixed internal selling price which is lower than the internal buying price. The profit margin between the internal selling and buying price is for prosumers renting the capacity of CES. Therefore, the profits of the ESCO consist of the NPV of energy revenue from both the internal trading profits in the planning horizon minus the initial CES deployment investment. \red{Many governments and municipalities offer incentives, subsidies, and/or favorable policies to promote energy efficiency and sustainability. These incentives can make ESCO projects financially attractive to investors, as they enhance the return on investment. }Thus this %cluster setting 
\red{business model} will attract more intermediary agents to invest in the energy service by setting the community microgrid and internal energy transaction business platform and thus will accelerate the transformation from the traditional electricity market to a decentralized market.

\subsection{\red{Summary of the proposed models}}
\red{In the following, we summarize the main properties of the proposed business models. The objectives are in line with the major incentives for the different decision makers. So, for example, in the IL model, the whole community acts as an investor and shares the overall investment and fuel consumption costs of the microturbines. In the CN model, the community acts as an investor and collects the revenue from the \red{electricity} sold, and pays for the imported \red{electricity}. Hence, the objective is to minimize the energy costs traded with the main grid. Finally, in the ESCO model, the objective is to maximize the profits of ESCO that acts as an investor and collects revenue from trading energy to the households. 

In the IL and the CN models, the whole community benefits as the internal energy produced within the community is shared among all households at no cost. On the other side, the incentive of the households to participate in the ESCO model is given by the lower internal buying price compared to the external real-time price thus resulting in a lower household's energy consumption cost.}
%TODO... justify further why these: optimal investment, min import (which might benefit the prosumers more than the consumers in the LEC. Is “exportation revenue” paid for charging the CES or export to the grid?

\begin{table}[H]
\centering \red{
\begin{tabular}{ll|lll}
\toprule
\multicolumn{2}{l|}{%Clustering 
Business Model}                               & Island    & Interconnected & ESCO   \\
\midrule
\multicolumn{2}{l|}{Objective}                                              & min investment          &  min imported    & max profit  \\
 & &\& fuel  costs         &   energy cost    &   \\[1ex]
\multicolumn{2}{l|}{Decision maker}                                              & community        & community    & ESCO \\[1ex]
\multicolumn{2}{l|}{Unmet demand}                                                        & not allowed        & main grid     & main grid  \\[1ex]  
\multicolumn{2}{l|}{Specific constraints}                                                        & meet full demand internally        & investment budget     & None \\ [1ex] 
\multicolumn{2}{l|}{Investment costs carried by}                                                        & community        & community     & ESCO \\  [1ex]
\multicolumn{2}{l|}{Energy consumption costs carried by}                                                        &community         & community      & households \\  
\bottomrule
\end{tabular}}
\caption{Comparison of the different properties of the three business models}
\label{tab:bat2}
\end{table}

\section{Mathematical model}\label{sec:formulation}
\subsection{Notations}

The community power network of a neighbourhood can be modeled as a bi-directed graph $G=(N,A)$, where $N=C\cup P$. The sets $C$ and $P$ represent the consumer households set and prosumer households set, indexed by $j\in C=\{1,2,...,m\}$ and $k\in P=\{1,2,\cdots,n\}$ respectively. The potential battery location set $B$
in a given  neighborhood is indexed by $i\in \{1,2,\cdots,b\}$ and the battery capacity set $L$ is indexed by $l\in\{1,2,\cdots,|L|\}$. 
\red{We assume that the community energy storage is deployed next to a prosumer's house, and hence  $B \subseteq P$. The indices and the corresponding sets used in our models are defined in Table \ref{tab:parameters}. }

\begin{table}[H]
\renewcommand{\tabcolsep}{0.1cm}
\begin{tabular}{ll}
\hline
Index/Set& Description                                                                                \\ \hline
$l$ & index for CES battery capacity, $l\in \{1,2,..., |L|\}$    \\
$i$  & index for battery location, $i\in \{1,2,...,b\}$\\
$j$ &index for consumer location, $j\in \{1,2,...,m\}$ \\
$k$ & index for prosumer location , $k\in \{1,2,...,n\}$ \\
$t$  & index for time period, $t\in T= \{1,2,...,|T|\}$       \\
$L$ & set of battery capacities     \\
$B$  & set of potential community  energy storage locations        \\
$C$ & set of consumer locations \\
$P$ & set of prosumer locations       \\
$T$  & set of time periods per day \\ \hline
\end{tabular}
\caption{Indicies and Sets used in the models.}\label{tab:sets}
\end{table}

\begin{table}[H]
\renewcommand{\tabcolsep}{0.1cm}
\begin{tabular}{ll}
\hline
Parameter& Description                                                                                \\ \hline
$\cc^l_i$ & cost of deployment of a CES at location $i \in B$ with capacity $l \in L$    \\
$\scc$  & installation/setup cost for a microturbine \\
$\fc$ & fuel cost  for microturbines  \\
$dist_{ij}$  & distance limit between consumer $j$ and CES $i$ of each cluster \\
$dist_{ik}$ & distance limit between prosumer $k$ and CES $i$ of each cluster \\
$d_j^t$ & electricity demand of consumer $j$ at time $t$ \\
$d_k^t$  & electricity demand of prosumer $k$ at time $t$ \\
$p_k^t$  & electricity generation of prosumer $k$ at time $t$ \\
$F_{max}$  & flow capacity limit of transmission lines between CES of node $i$ and consumer node $j$ \\
$P^{l}_{dis}$  & nominal discharging power of a CES with capacity $l$ \\
$P^{l}_{ch}$     &  nominal charging power of a CES with capacity $l$ \\
$Pmt_{max}$   & maximum power generation of a microturbine            \\
$S^l_{min}$ &	minimum energy charged in a CES with capacity $l$ \\
$S^l_{max}$	& maximum energy charged in a CES with capacity $l$   \\       
$\eta_d$ &	discharging efficiencies of CES  \\
$\eta_c$	& charging efficiencies of CES  \\  
$\Delta t  $   &length of each time step $t$ \\
$\pi_{buy}$	&cluster internal electricity buying price per unit for household\\
$\pi_{sell}$&cluster internal electricity selling price per unit for prosumer\\
$\pi_{ex}$	&price per unit of selling electricity to main grid in the interconnected clustering model\\
$\pi_{im}^t$  	&price per unit of buying electricity at time $t$ for households at real-time market\\
%$\pi_{wlt}^t$	&price per unit of buying electricity at time $t$ for ESCO at wholesale market\\
%$\mu$& profit margin of ESCO \\
%$n$ & number of CES long-term planning years\\
$r$ & discount rate of NPV\\
$\rho$ &  saving factor of electricity bills to invest in CESs \\
$\alpha$ & NPV factor of the planning years \\
\hline
\end{tabular}
\caption{Parameters used in the Model}\label{tab:parameters}
\end{table}

\red{In addition, potential microturbine locations are also determined by the set $B$. In this case, the same network infrastructure can be used to transport energy from the battery/microturbine site to the households. Therefore, it is reasonable to assume that microturbines and CES share the same potential locations.}
Moreover, there are $|L|$ different types of CES that can be deployed, each with a given charging/discharging capacity $P^l_{ch}$ and $P^l_{dis}$, respectively, $l \in L$. Installing a CES of type $l$ at location $i$ incurs a cost of ${\cc}^l_i$.

The household's load and prosumer's generation profile data are obtained from smart meter measurements at each time-slot. For every $t \in T$, we are given the \red{following data defined in Table \ref{tab:parameters}: $d_j^t$ , $d_k^t$, and $p_k^t$.}  %[-3ex]
%\begin{itemize}
%    \item $d_j^t$  electricity demand of consumer $j \in C$ at time $t$, \\[-4ex]
%    \item $d_k^t$   electricity demand of prosumer $k \in P$ at time $t$, \\[-4ex]
%    \item  $p_k^t$   electricity generation of prosumer $k \in P$ at time $t$. 
%\end{itemize}
%A consumer's $j \in C$ electricity power demand is 
\red{The power demand of consumer $j \in C$ is} represented by a load profile determined by the vector $d$ above.
To set the prosumer power flow operation mode (charging/discharging), we use the following function to judge if the prosumer generates electricity or consumes electricity:
\begin{align*}
    [\delta]^{+} & =\left\{
                \begin{array}{ll}
                  \delta \quad \delta \ge 0\\
                  0 \quad\text{otherwise.}
                \end{array}
              \right.
\end{align*}
Hence, a prosumer's $k$ net demand at time $t$ is determined by the difference between their load and the PV generation profile, given as $[d_k^t-p_k^t]^+$. Similarly the prosumer $k$ net charging at time $t$ is determined as $[p_k^t-d_k^t]^+$.

%I checked from the code it is every fifth prosumer.}}   The problem can be formulated as a capacitated facility location problem where the decision variables are defined as follows:

Our models make use of the following network design variables, that determine the underlying infrastructure, location of CES and their configuration:
\begin{align*}
    y_i^l & =\left\{
                \begin{array}{ll}
                  1 \quad\text{if a CES of capacity $l$ is deployed at location $i$,\quad $i\in B,l\in L$}\\
                  0 \quad\text{otherwise.}
                \end{array}
              \right.           \\
                  %\textcolor{red}{y_{i,mt}} 
                  {\xi_{i}} & =\left\{
                \begin{array}{ll}
                  1 \quad\text{if a microturbine is deployed at location $i$,\quad $i\in B$}\\
                  0 \quad\text{otherwise.}
                \end{array}
              \right.           \\
    x_{ij} & =\left\{
                \begin{array}{ll}
                  1 \quad\text{if the CES at node $i$ is allocated to consumer $j$,\quad $i\in B,j\in C$}\\
                  0 \quad\text{otherwise.}
                \end{array}
              \right.\\
    z_{ik} & =\left\{
                \begin{array}{ll}
                  1 \quad\text{if the CES at node $i$ is allocated to prosumer $k$,\quad $i\in B,k\in P$}\\
                  0 \quad\text{otherwise.}
                \end{array}
              \right.
\end{align*}

\red{Table \ref{tab:variables} lists all continuous variables used in our model.} They represent the DC network power flow between batteries and households (which is common for all three models). \red{Note that we assume line losses are negligible and can be ignored when considering short distances on the low voltage network.} For the island model (IL), we also have power flow variables between microturibines and pro/consumers, and for the interconnected (CN) and ESCO model, we have power flow variables modeling the energy exchange between households and the main grid \citep{CALVILLO2015}. For prosumers, the flow direction is supposed to be bi-directional because of its double purpose: acting as both an energy consumer and a generator. 

\begin{table}[H]
\renewcommand{\tabcolsep}{0.1cm}
\begin{tabular}{lll}
\hline
Vars   &Model        & Description                                                                                \\ \hline
$f_{ij}^{t}$ & all & power flow from CES of node $i$ to consumer $j$ at time $t$; $i \in B, j \in C, t \in T$\\
$g_{ik}^{t}$ & all & power flow from CES of node $i$ to prosumer $k$  at time $t$; $i \in B, k \in P, t \in T$ \\
$g_{ki}^{t}$ & all & power flow from prosumer $k$ to CES of node $i$ at time $t$;  $k \in P, i \in B, t \in T$ \\
$S_i^t$		& all & the amount of electricity charged into the CES of node $i$ at time $t$; $i \in B, t \in T$ \\
$u_{ij}^t$	& IL & power flow from microturbine of node $i$ to consumer $j$ at time $t$; $i \in B, j \in C,  t \in T$\\
$u_{ik}^t$	& IL & power flow from microturbine of node $i$ to prosumer $k$ at time $t$; $i \in B, k \in P, t \in T$\\
$u_{j}^t$    & CN/ESCO & electricity that consumer $j$ imports from main grid at time $t$; $j \in C, t \in T$\\
$u_{k}^{t+}$	& CN/ESCO & electricity that prosumer $k$ imports from main grid at time $t$; $k \in P, t \in T$\\
$u_{k}^{t-}$	& CN/ESCO & electricity that prosumer $k$ exports to main grid at time $t$;  $k \in P, t \in T$\\ 
 \hline
\end{tabular}
\caption{Power flow variables \label{tab:variables}}
\end{table}
\noindent Note that $S_i^t$ is the amount of electricity stored in a battery, which is equal to the nominal capacity of the battery times its State of Charge (SOC). This variable is acting as the inventory level of a warehouse in a distribution network, which is inter-temporally connected because of charging and discharging operations.
\subsection{Island model}	
In the island model, it is crucial for utility companies to dispatch the electricity that prosumers generate at the microgrid level via an EMS and \red{in the case of a shortfall to supply} households with electricity from microturbines to guarantee the self-sufficiency of each microgrid. 

The objective function for the community power network design problem of the Island model (IL) can be formulated as:
%the following mixed integer programming model:
\begin{equation}
\min \quad \sum_{i\in B}{\sum_{l\in L}\cc_i^ly_i^l}+\scc \sum_{i \in B} \xi_{i} +  {\alpha} \cdot \fc \sum_{t\in T}\sum_{i\in B}Pmt_{i}^t,   \label{of1}
\end{equation}
where $\xi_{i}$ indicates whether a microturbine is installed at location $i$ or not, $y^l_i$ indicates whether a CES of type $l$ is installed at $i$ and $Pmt_{i}^t$ is the amount of energy produced by the microturbine located at node $i$ at period $t$. 
%In this formula, the first term corresponds to CES installation cost, and the total cost function of microturbines is: $\scc \sum_{i \in B} \xi_{i} +\fc \sum_{t\in T}\sum_{i\in B}P_{i}^t$ where $\scc$ and $\fc$ are the installation/setup and fuel cost coefficients respectively. 
In order to choose the optimal capacity and location from the potential battery deployment nodes, the objective function \eqref{of1} minimizes the investment cost of CES (represented by the first term in the objective function), as well as the installation and the discounted fuel cost of the microturbine operations, given by the second and the third term, respectively. Hereby, $\scc$ and $\fc$ are the installation/setup and fuel cost coefficients for the microturbines, respectively, and {$\alpha$} is the NPV factor of the planning years given as  %\blue{$\alpha=\sum_{n=1}^{\#years}\frac{365}{(1+r)^n}$} %where 365 is the total days in a year, 
{$\alpha=\sum_{n=1}^{\#years}\frac{8760/|T|}{(1+r)^n}$}
to adjust the daily planning cost from the objective function.\footnote{{We assume that households' load and prosumers' generation profile data are given on an hourly basis.}} The overall microturbine generation power at each timestamp $t$ is  given in \eqref{of2}. The power flow $u$ can be sent out from location $i \in B$, only if a microturbine has been installed at this location, while the maximum power generation per microturbine is limited to 
%$P_{mt,max}$. 
$Pmt_{max}$,  see \eqref{MT}:
\begin{subequations}
\label{subeq:mt}
\begin{align}
Pmt_{i}^t & = \sum_{j\in C}u_{ij}^t+\sum_{k\in P}u_{ik}^{t} && \quad t\in T, i \in B \label{of2}\\
\sum_{j\in C}u_{ij}^t+\sum_{k\in P}u_{ik}^t & \le %P_{mt,max}y_{i,mt}
Pmt_{max} \xi_i
%\sum_{l\in L} y_i^l, 
 && \quad t\in T,i\in B  \label{MT}
\end{align}
\end{subequations}
%The initial deployment costs and fuel costs are paid by utility companies while accommodating the constraints that are applicable to the system as follows.
There are three groups of constraints that constitute our MILP model: those describing the infrastructure (network-design constraints), those guaranteeing flow-balance in the network (network power flow constraints) and those guaranteeing that the capacity and charging/discharging at CES are respected (CES power flow  constraints). 
Network-design constraints for the IL model are given as follows:
\begin{subequations}
\label{subeq:NetDes}
\begin{align}
\sum_{i\in B} x_{ij}  & \leq 1\ && \quad j\in C \label{eq:x}\\
\sum_{i\in B} z_{ik} & \leq 1\ && \quad k\in P  \label{eq:z}\\
\sum_{l\in L}{y}_i^l & \leq 1 &&  \quad   i\in B \label{eq:yl}\\
x_{ij} & \le \sum_{l\in L}y^l_i && \quad j\in C, \quad  i\in B \label{eq:xy}\\
z_{ik} & \le \sum_{l\in L}y^l_i && \quad k\in P, \quad  i\in B \label{eq:zy}\\
x_{ij} & =0 && \quad   i\in B, j \in C:  dist_{ij} > Dist_{max}  \label{eq:x0}\\
z_{ik} & =0 && \quad  i\in B, k \in P:  dist_{ik} > Dist_{max}  \label{eq:z0}\\
y_i^l & \in\{0,1\} && \quad  i\in B,\ l\in L 	\\
x_{ij} & \in\{0,1\} && \quad i\in B,\ j\in C\\
z_{ik} & \in\{0,1\} && \quad  i\in B,\ k\in P \label{integer}\\
\xi_i & \in\{0,1\} && \quad  i\in B \label{integer_xi}
\end{align}
\end{subequations}

Constraints \eqref{eq:x} and \eqref{eq:z} ensure that each household is connected to at most one battery. Constraint \eqref{eq:yl} \red{selects} the capacity $l$ of each battery, and if there is no battery deployed at node $i$, $\sum_{l \in L} y_i^l$ is set to 0. \red{Constraints \eqref{eq:xy} and \eqref{eq:zy} ensure that if a battery %with capacity $i$ 
is deployed at node $i$ then a capacity for that battery needs to be selected}. Additionally, $dist_{ij}$ and $dist_{ik}$ are the distances between a battery and a consumer/prosumer respectively and $Dist_{max}$ is the maximum allowable distance between a battery and a household in one cluster. The options for clustering will be limited by the distance constraints \eqref{eq:x0} and \eqref{eq:z0}.

We now present the network power flow constraints. In order to ensure the balance between power supply and demand for consumers in each cluster, the  power conservation constraint \eqref{eq:d}
states that in each time period $t$, the total demand of a consumer $j \in C$ (given as $d_j^t$) is satisfied using flow that emanates from microturbines or from CES. 
Similarly, the power demand and supply conservation constraints for prosumers are given by \eqref{eq:pd}-\eqref{eq:dp}.
Constraint \eqref{eq:Pdisf} restricts the power flow capacity from each location $i \in B$ to each consumer $j \in C$, same as the constraint \eqref{eq:Pdischg} for prosumers while charging CES and discharging CES/microturbine.
$F_{max}$ represents the flow capacity of transmission lines between a node $i \in B$ and a household node. For given $P^{l}_{dis},P^{l}_{ch}$ (nominal discharging and charging power limits of a CES with capacity $l$), constraints \eqref{eq:Pdis} and \eqref{eq:Pch} impose battery power configuration restrictions on network. 
\begin{subequations}
\label{subeq:power}
\begin{align}
\sum_{{i}\in B} (u_{ij}^t+f_{ij}^t) & =d_j^t  && \quad   t\in T, j\in C  \label{eq:d}\\
\sum_{{i}\in B}(u_{ik}^{t}+ g_{ik}^{t}) & =[d_k^t-p_k^t]^+ && \quad  t\in T, k\in P  \label{eq:pd}\\
\sum_{{i}\in B} g_{ki}^{t} & \le [p_k^t-d_k^t]^+ && \quad  t\in T, k\in P  \label{eq:dp} \\
%\end{align}
%\end{subequations}
%In constraint \eqref{eq:d}, $d_j^t$ refers to the electricity demand of consumer $j$ at time-slot $t$, similarly, $d_k^t$ refers to the electricity demand of prosumer $k$ at time $t$ and $p_k^t$ is the electric energy produced by the prosumer $k$ at time $t$. 
%The other network power flow constraints for CES and microturbine configuration are as follows:
%\begin{subequations}
%\label{subeq:powerflow}
%\begin{align}
f_{ij}^t+u_{ij}^t & \le F_{max}\ x_{ij}  && \quad t\in T,i\in B,\ j\in C \label{eq:Pdisf}\\
g_{ik}^{t}+g_{ki}^{t}+u_{ik}^{t} & \le F_{max}\ z_{ik} &&  \quad  t\in T,i\in B,\ k\in P \label{eq:Pdischg}\\
\sum_{j\in C}{f_{ij}^t}+\sum_{k\in P} g_{ik}^{t} & \le \ \sum_{l\in L} P^{l}_{dis} y_i^l && \quad  t\in T,i\in B \label{eq:Pdis}\\
\sum_{k\in P} g_{ki}^{t} & \le \ \sum_{l\in L} P^{l}_{ch} y_i^l && \quad  t\in T,i\in B  \label{eq:Pch}\\
f_{ij}^t,u_{ij}^t & \geq0 &&  \quad t\in T,i\in B,\ j\in C  \\  	
g_{ik}^{t},g_{ki}^{t},u_{ik}^{t} & \geq0 &&  \quad  t\in T,i\in B, k\in P \label{continuous}
\end{align}
\end{subequations}

Finally, we present the CES power flow constraints. Since each CES plays a role of energy warehouse that coordinates the electricity charging and discharging planning, the inter-temporal and physical constraints for the CES are as follows (see also \cite{chen2011sizing}):

\begin{subequations}
\label{subeq:CES}
\begin{align}
S_i^{t} & =S_i^{t-1} + \eta_c  \Delta t   \sum_{k\in P} g_{ki}^{t} -\frac{ \ \Delta t}{\eta_d } (\sum_{k\in P} g_{ik}^{t}+\sum_{j\in C}{f_{ij}^{t}}) && \quad t\in T ,  i\in B 	 \label{eq:SOC} \\
\sum_{l\in L} S_{min}^{l} y_i^l  & \le S_i^t  \le \sum_{l\in L}  S_{max}^{l} y_i^l && \quad  t\in T,i\in B   \label{eq:Sm}\\
S_i^0 & = S_i^{24} && \quad i\in B \label{eq:SOC0} \\ 
%&\bissan{\text{what about the next day and the day after, how is this handled in the code?}} \notag\\
S_i^t & \geq0 && \quad  t\in T,i\in B\label{S}
\end{align}
\end{subequations}

%$\omega$ is the initial SOC proportion of CES capacity. 
The parameters $\eta_d$ and  $\eta_c$ are the discharging and charging efficiencies for a CES deployed at node $i$.
$\Delta t $ is the length of the time step $t$. Constraint \eqref{eq:SOC} shows the intertemporal state of charge connected with charging and discharging flows. CES capacity restriction is given by constraint \eqref{eq:Sm}. In this work, we assume that the charging and discharging process of each CES in one day is a cycle, which is given in constraint \eqref{eq:SOC0}.

\subsection{Interconnected model}
In the interconnected model (CN), communities are connected to the main grid so that the utility company can provide electricity when the prosumers of one community cannot supply enough energy for all households. At the same time, if the community generates excessive energy and the CES are charged \red{to maximum capacity,} the prosumers can also export electricity to the main grid. 
%\red{The total electricity bill savings of one community are shared by the households in the community according to the electricity they consume, and the ownership of CES is all the households connected to the network.}\footnote{Do we need this sentence, we may be attacked by the reviewers because of this assumption, unless we have some reference to cite where this is indeed the case? \bissan{I agree}}  Therefore 
The objective function in this interconnected model is to maximize the NPV revenue that a community obtains from exporting energy to the main grid minus the expenses that they spend on importing energy from the main grid. Nowadays because of the insufficient amount of solar panel installation of prosumers, the energy sharing communities cannot be net energy generators for the main grid. Thus we modify the objective function to minimize the energy cost as follows:
\begin{align}
min \quad {\alpha} \cdot \sum_{t\in T}\left[\left(\sum_{j\in C}u_j^{t}+\sum_{k\in P}u_k^{t+}\right)\pi_{im}^t- \sum_{k\in P} u_k^{t-}\pi_{ex}\right],
\end{align}
\red{where {$\alpha$} is the NPV factor of the planning years given as in Section \ref{sec:IL} to adjust the daily planning cost from the objective function.} The parameter $\pi_{im}^t$ is the price per unit of buying electricity from the main grid in real-time electricity market at time $t$ and $\pi_{ex}$ is the {fixed} price per unit of selling electricity from prosumers to the main grid in the electricity market.
%at time $t$. 
When the generation of a \red{community cannot serve its own demand,} it buys energy from the main grid. Otherwise it can export and sell it to the utility company. In this setting, we assume that the importing price is time-related to incentivize the households to shift their electricity consumption from high price to low price hours, in order to reduce their expenditures and to lead the least efficient power plants in the main grid to stop production. {The fixed exporting price is relatively lower than the importing price.} This pricing mechanism aims at improving households' self-consumption and self-sufficiency through demand side management and coordination of EMS.

Let $D$ represent the electricity bill cost  over the given planning horizon when there is no CES deployed. This cost consists of the energy expenses for both consumers and prosumers, which is the amount of money they pay the utility company for the energy they import from the main grid. The prosumers can also sell the surplus to the main grid after providing themselves with the energy they need to get some profits from the utility company, as they serve as the electricity generators in the main grid for the utility company. Hence, this constant $D$ is calculated as: \begin{equation*}
D={\alpha} \cdot \sum_{t\in T}\left({\sum_{j\in C}d_j^{t}\pi_{im}^t }
+{\sum_{k\in P}{[d_k^{t}-p_k^{t}]^+}\pi_{im}^t}
-{\sum_{k\in P}{[p_k^{t}-d_k^{t}]^+}{\pi_{ex}}}\right).
\end{equation*}
In our model, in order to incentivize the community households to invest in CES, we set the available CES investment budget as a certain percentage saving factor $\rho$ with respect to $D$ (recall that  $cc_i^l$ represents the cost of deployment of a battery with capacity $l$):
\begin{equation}
\sum_{i\in B}{\sum_{l\in L}\cc_i^ly_i^l}  \leq \rho D\label{budget}
\end{equation}

%The available investment budget of the storage equipments for one neighborhood is given as $\rho D$ in constraint \eqref{budget}, where saving factor $\rho$ is given by households: decision maker in this clustering model and $cc_i^l$ represents the cost of deployment of a battery with capacity $l$.

Also in this model we distinguish between: network-design constraints, network power flow constraints and CES power flow constraints. Network design decisions are modeled using binary variables $(x,y,z)$ introduced above, and the associated network-design constraints are given by \eqref{eq:x}-\eqref{integer}. The CES power flow constraints \eqref{subeq:CES} remain the same. Since the communities are connected to the main grid, the network power flow constraints are modified as follows:
\begin{subequations}
\label{subeq:inter}
\begin{align}
u_j^t+\sum_{{i}\in B} f_{ij}^t & =d_j^t && \quad   t\in T, j\in C \label{eq:d1}\\
u_{k}^{t+}+\sum_{{i}\in B}g_{ik}^{t} & =[d_k^t-p_k^t]^+ && \quad  t\in T, k\in P  \label{eq:pd1}\\
u_{k}^{t-}+\sum_{{i}\in B} g_{ki}^{t} & = [p_k^t-d_k^t]^+ && \quad  t\in T, k\in P  \label{eq:dp1} \\
u_{j}^t & \le F_{max} &&  \quad t\in T,\ j\in C \label{eq:uUp} \\
u_{k}^{t+}+u_{k}^{t-} & \le F_{max} &&  \quad t\in T,\ k\in P\\
f_{ij}^t & \le F_{max}\ x_{ij}  && \quad t\in T,i\in B,\ j\in C \label{eq:Pdisf1}\\
g_{ik}^{t}+g_{ki}^{t} & \le F_{max}\ z_{ik} &&  \quad  t\in T,i\in B,\ k\in P \label{eq:Pdischg1}\\
u_{k}^{t+}, u_{k}^{t-}, u_{j}^{t} & \geq0 &&  \quad  t\in T,i\in B, k\in P  \label{eq:urange}\\
& \text{ \eqref{eq:Pdis}-\eqref{continuous}} && \notag
\end{align}
\end{subequations}
Constraints \eqref{eq:d1}-\eqref{eq:pd1} guarantee that each household can satisfy its energy demand from the main grid as well (see variables $u_j^t$ and $u_k^t$). The balance constraints \eqref{eq:dp1} {determine the net energy production of prosumers as well as the energy export to the grid.} %\red{are also where the demand of consumers and the net energy production of prosumers are decided on both community network and wider network power flow.}\footnote{Are we happy with this sentence?}
Constraints \eqref{eq:uUp}-\eqref{eq:urange} provide upper and lower bounds on the power flow variables. 

%The other constraints remain the same as for the island model, given as \eqref{eq:x}-\eqref{integer}, \eqref{eq:Pdis}-\eqref{S}.

\subsection{ESCO model}

To help households reduce their electricity bills by stocking the energy that prosumers produce and then dispatch to each household, in this setup the ESCOs charge households for the electricity they use from the battery at a fixed internal buying price ($\pi_{buy}$) and pay prosumers for the energy they produce ($\pi_{sell}$). {From the perspective of ESCO,} the objective function is thus to maximize the NPV of the internal energy trading profits minus the investment cost for shared batteries:
\begin{flalign}
max\quad {\alpha} \cdot \sum_{t\in T}\sum_{i\in B}\left[ \left( \sum_{j\in C}f_{ij}^t+\sum_{k\in P}g_{ik}^t \right)\pi_{buy}-
\sum_{k\in P}g_{ki}^t\pi_{sell} \right] 
%+\mu((\sum_{j\in C}u_{j}^t+\sum_{k\in P}u_{k}^{t+})\pi_{im}^t+\sum_{k\in P}u_{k}^{t-}\pi_{ex}^t)] \nonumber 
-\sum_{i\in B}{\sum_{l\in L}\cc_i^ly_i^l}
\end{flalign}
\red{where {$\alpha$} is the NPV factor of the planning years defined in Section \ref{sec:IL}.} The parameter $\pi_{buy}$ is the internal energy buying price per unit from discharging the battery, and $\pi_{sell}$ is the internal selling price for prosumers charging the CES. The objective is to maximize the NPV of ESCO's profits minus the initial investment for CES. The network-design constraints are   \eqref{eq:x}-\eqref{integer}, the network power flow constraints are  \eqref{eq:Pdis}-\eqref{continuous},  \eqref{subeq:inter}. 
Finally, the CES power flow constraints are \eqref{subeq:CES}.

\subsection{Benders Decomposition}
The three MILP problems described in the sections above, exhibit a structure that can be exploited using a Benders decomposition approach. The goal is to project out all continuous variables and keep only the binary variables in the master problem. 

\subsubsection*{Decomposing the Island Model} Starting from the island model, the problem can be reformulated as follows:
\begin{align*}
%(\text{MP}^{IL}) \qquad 
\min & \quad \sum_{i\in B}{\sum_{l\in L}\cc_i^ly_i^l}+\scc \sum_{i \in B} \xi_{i} + \Theta^{IL}   \\
& \text{s.t. $(x,y,z,\xi)$ satisfy \eqref{subeq:NetDes} }  \\
& \Theta^{IL} \ge \Phi^{IL}(x,y,z,\xi)
% (\text{MP}^{CN}) \qquad \min & \quad \Phi^{CN}(x,y,z)  \qquad \\
% & \text{s.t. $(x,y,z)$ satisfy \eqref{eq:x}-\eqref{integer}, \eqref{eq:Sm}-\eqref{S} and \eqref{budget}} \\
% (\text{MP}^{ESCO}) \qquad \max & \quad \Phi^{ESCO}(x,y,z) -\sum_{i\in B}{\sum_{l\in L}cc_i^ly_i^l}   \\
% & \text{s.t. $(x,y,z)$ satisfy \eqref{eq:x}-\eqref{integer} and \eqref{eq:Sm}-\eqref{S}} \\
\end{align*}
where for any given allocation of households to batteries (determined by the vector $(x,y,z,\xi)$), the function $\Phi^{IL}(.)$ calculates the
\red{total fuel cost} to satisfy energy demand in the island model. If the vector $(x,y,z,\xi)$ results into an infeasible configuration with respect to SOCs and/or customer demands, we will assume that $\Phi^{IL}(.)= \infty$. The auxiliary variable $\Theta^{IL}$ is bounded from below by these costs, and at the optimum, $\Theta^{IL}$ is equal to the the fuel costs for the given network design decision:
\begin{align*}
\Phi^{IL}(x,y,z,\xi)= & \min_{(f,g,u,S) \ge 0} {\alpha} \cdot \sum_{t\in T}\sum_{i\in B} \fc (\sum_{j\in C}u_{ij}^t+\sum_{k\in P}u_{ik}^{t})  \\
& \text{s.t. $(f,g,u,S) $ satisfy \eqref{subeq:mt}, \eqref{subeq:power}-\eqref{subeq:CES}}\\
\end{align*}

%  or the profit (or deficit) obtained from selling the suprplus energy minus the import of energy from the grid for the whole time period $T$ for the other two models. 
Because of the SOC constraint \eqref{eq:SOC} that links the time-periods $t-1$ and $t$, the associated Benders subproblem cannot be decomposed per time period, and thus it represents one large LP.  
Hence, only a single cut can be generated in each iteration after solving the Benders subproblem. The result of this LP will induce either a Benders feasibility cut (if the solution of the master problem determined by the vector $(x,y,z,\xi)$ renders the subproblem infeasible), or a Benders optimality cut (if the value $\Theta^{IL}$ of the master problem is smaller than the solution value of the subproblem's LP).
For more details concerning implementation of Benders decomposition, see \cite{cplex-benders}. In our implementation, to test the performance of the Benders decomposition method, we use CPLEX automatic Benders decomposition algorithm which automatically associates all binary variables to the master problem and all other continuous variables to the subproblem. The algorithm uses some of state-of-the-art stabilization and acceleration techniques, such as cut-loop separation and in-out approach, see, e.g. \cite{Ben-AmeurN07,FischettiLS16, FischettiLS17}. In what follows, we explain how the decomposition is guided by the value function reformulation for the remaining two models. 

\subsubsection*{Decomposing the Interconnected Model}
Similarly, for the interconnected model, we reformulate the problem as
\begin{align*}
\min & \quad \Theta^{CN} \qquad \\
 & \text{s.t. $(x,y,z)$ satisfy \eqref{eq:x}-\eqref{integer},  \eqref{budget}} \\
 & \Theta^{CN} \ge \Phi^{CN}(x,y,z)
\end{align*}
where the value of the function $\Phi^{CN}(x,y,z)$ is obtained by solving the following LP:
\begin{align*}
\Phi^{CN}(x,y,z)= &
\min \quad  {\alpha} \cdot \sum_{t\in T}\left[\left(\sum_{j\in C}u_j^{t}+\sum_{k\in P}u_k^{t+}\right)\pi_{im}^t- \sum_{k\in P} u_k^{t-}\pi_{ex}\right] \\
& \text{s.t. $(f,g,u,u^+,S) $ satisfy  \eqref{eq:Pdis}-\eqref{continuous}, \eqref{subeq:CES}}, \eqref{subeq:inter}\\
\end{align*}

\subsubsection*{Decomposing the ESCO Model}
Finally, the decomposition of the ESCO model is guided by the following reformulation:
\begin{align*}
\max & \quad \Theta^{ESCO} - \sum_{i\in B}{\sum_{l\in L}\cc_i^ly_i^l}\qquad \\
 & \text{s.t. $(x,y,z)$ satisfy \eqref{eq:x}-\eqref{integer}}
 %,  \eqref{budget}} 
 \\
 & \Theta^{ESCO} \le \Phi^{ESCO}(x,y,z)
\end{align*}
where the value of the function $\Phi^{ESCO}(x,y,z)$ is obtained by solving the following LP:
\begin{align*}
\Phi^{ESCO}(x,y,z)= &
\max \quad   
{\alpha} \cdot \sum_{t\in T}\sum_{i\in B}\left[ \left( \sum_{j\in C}f_{ij}^t+\sum_{k\in P}g_{ik}^t \right)\pi_{buy}-
\sum_{k\in P}g_{ki}^t\pi_{sell} \right] 
\\
& \text{s.t. $(f,g,u,u^+,S) $ satisfy  \eqref{eq:Pdis}-\eqref{continuous}, \eqref{subeq:CES}, \eqref{subeq:inter}} \\
\end{align*}

\section{Case study: Cambridge, MA}\label{sec:usecase}
In this section, we present and analyze our results for the three proposed \red{business models.} These results are the numerical solutions for three different stakeholders who minimize their costs or maximize their profits: utility companies, prosumer groups, and Energy Service Companies, respectively. As explained in Section \ref{sec:3models}, these three entities have different %clustering 
optimization goals which are formulated in Section \ref{sec:formulation}. We apply these formulations to a case study and then compare the clustering and CES deployment results. 

We use a real case study with the geographical data as well as the energy consumption and generation profiles taken from \citep{barbour2018community}, which represents an electricity network of Cambridge Massachusetts, an area that covers 4574 households of consumers and prosumers. To obtain representative 
%days' 
{daily} data for both energy demand and generation, we converted the 15-min resolution data of each day into hourly resolution data. For each household used, the corresponding demand is summed over each group of 4 consecutive 15-minute demand values. Using this hourly consumption, the mean demand per hour of the day, across the total of 30 days, is then taken, yielding 24 values per household, representing the hourly electricity consumption on an average day {of the given period of 30 days}. {We note that the Cambridge demand distributions are similar
throughout the year (as shown in \cite{barbour2018community}).
{Following the study of \cite{barbour2018community} \red{where they show that the demand variation over the year is limited}, we also restrict the demand data to one representative month (July in this case).} At the end of this section, we discuss the effect of varying the demand data on our estimated operational costs.} 
%\footnote{Didn't we look separately at summer vs winter consumption/production?} 

\blue{In this section, we analyze the results of the case study using a small neighborhood of 120 households consisting of 40 prosumers and 80 consumers. The set of prosumers and consumers are sampled from the given 4574 households in a way that prevents generating neighborhoods where the households are located far-off, in particular for small instances.} \red{In practice, it is reasonable to divide a given geographical area into different communities that are geographically close and to optimize separately over these communities as they can be considered independent. It is also not possible to connect households with batteries from distant communities due to the physical limitations of the transmission lines.} \blue{Additionally, from the dataset, we know the households' consumption and production levels. Thus we choose consumers from the ones whose production level is zero and prosumers from the ones with a non-zero production level.} The exact topology of the power network is not openly available due to security concerns, but the geographical road network which connects all the households in one neighborhood can provide an approximation of a real community power network. {As in \citep{barbour2018community},} we assume that the geographical road network is the approximation of the real electrical distribution network, thus we use the open source routing machine, OSRM \cite{OSRM}, to calculate the shortest path in the road network using the real geographical coordinates and to represent the electricity network cable length between two households.

\begin{table}[H]
\centering
{
\begin{tabular}{@{}lllll@{}}
 \toprule       
& HES&\multicolumn{3}{l}{CES Configurations} \\  \midrule
Battery Size (kWh)    & 13.5  & 75     & 150   & 250  \\ 
Price(\$)            & 15,175 & 33,000  &  65,250  &107,500  \\
Battery Power(kW)    &5.5 &500 &500 &500 \\
\midrule
& \multicolumn{4}{l}{Other Battery Parameters} \\  \midrule
Fuel Price(\$/kWh)   &\multicolumn{4}{c}{0.5} \\
Life Cycle(years)  & \multicolumn{4}{c}{10} \\
$SOC_{min}$  & \multicolumn{4}{c}{10\%} \\
$SOC_{max}$ & \multicolumn{4}{c}{85\%} \\
$\eta_d=\eta_c$ & \multicolumn{4}{c}{95\%} \\ 
\bottomrule                         
\end{tabular}}
\caption{Overview of battery configurations}
\label{tab:bat}
\end{table}

The characteristics and parameter values for HES and CES are presented in Table \ref{tab:bat}.%\footnote{\bissan{ eCamion it still exists but the specs are no longer available and I could not find a paper that reference them (we can refer to my student's thesis that used these values https://uwspace.uwaterloo.ca/handle/10012/14628 or just state we have these values from 2019, what do you think?} Yes, let's cite your student and mention these are prices from 2019, along with the energy prices. } 
We adopt the HES settings from Sonnen Company \cite{Sonnen,SonnenPrice} which is used for single household energy storage setting and the CES settings are from a Canadian CES production company ``eCamion'' \cite{eCamion}. Three different capacities are considered for CES as shown in Table \ref{tab:bat} \cite{affleckthesis}. The optimal lifespan of lithium-ion batteries is set to 10 years \cite{2017life}. The potential CES locations set $B$ is a subset of $P$, and we assume that the number of the potential CES deployment locations is 20\% of the total number of prosumers.

\subsection{Data Generation}\label{sec:data_gen}
\red{The choice of 40 prosumers and 80 consumers samples from the original dataset was taken such that the households were close in terms of proximity and thus their geographical location was taken into account.} The identifiers that represent the households' coordinates data are then assigned to prosumers, consumers, and potential battery locations within the district based on the ratio of prosumers to consumers. {Because of the solar energy generation deficiency, we assume the ratio of prosumers to consumers in the neighborhood is either 2 or 3. %Thus, prosumer and consumer identifiers are assigned to the neighborhood household coordinates according to this ratio. 
%(i.e.,. 1/1, 2/1, 3/1, etc.). 
In this way of coordinates assignment, we can get a pre-defined network with evenly allocated prosumer and consumer locations in each neighborhood, which is essential in our model settings to avoid generating neighborhoods where there are only consumers or prosumers.

To generate the consumers' and prosumers' monthly load and generation profiles, we use the monthly electricity data from July (summer month) and December (winter month), 2015 in Cambridge, MA \citep{barbour2018community}. This provides a source of 15-min resolution electricity load profile. Two sets of solar generation profiles with the same temporal resolution for prosumers are also provided for summer days and winter days. A daily load profile and PV generation (winter and summer) data of one arbitrary prosumer is shown in Figure \ref{fig:gem}.  
\begin{figure}[H]\begin{center}
     \centering
\includegraphics[width=0.8\textwidth]{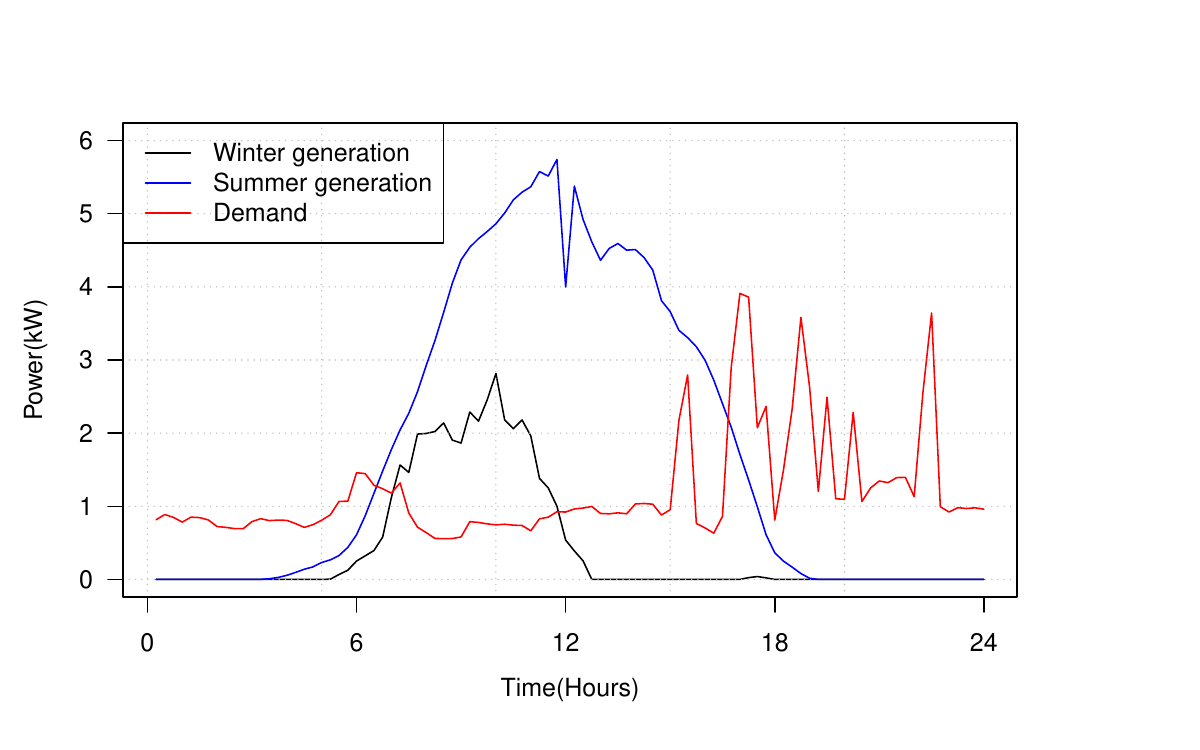}
\caption{One arbitrary prosumer's daily demand and power generation.}
\label{fig:gem}
\end{center}
\end{figure}

 %and we  get the average of one month. 
As shown in Figure \ref{fig:gem}, the winter generation is not sufficient compared to the load profile, thus we consider summer generation data in our case study. To maximize the CES utilisation efficiency, we assume that the energy conversion efficiency can be raised with the future advanced photovoltaic technology, therefore, we doubled the generation profile of the prosumers. %The households' electricity profiles data are then assigned to a selected district in Cambridge, MA, and the geographical data generation process is described in the Appendix I.
Our real-time electricity prices $\pi_{im}^t$ in Massachusetts is adopted from ISO New England \cite{Pricing} for the day of June 18th, 2019. We assume that the internal trading price set by ESCO is $\pi_{buy}=0.19224$ \$/kWh, which is the lowest importing external energy price. It is widely understood that at current US prices neither batteries nor PV are economic without subsidies \citep{2013prospects}. The rapid diffusion of solar PV has already resulted in operability issues and grid disruption in numerous markets \cite{parag2016electricity}. To relieve the operational burden of the main grid, we set the internal buying price to be lower than external real-time electricity price to make the internal transaction more competitive compared to external energy from the main grid, while the selling price for prosumers in internal exchange is the same as exporting to main grid price $\pi_{sell}=\pi_{ex}=0.05 $ \$/kWh. This price is also consistent with the export price in interconnected model. Figure \ref{fig:price} shows the electricity prices we adopt in our models.
\begin{figure}[!h]
\begin{center}
\includegraphics[width=0.8\textwidth]{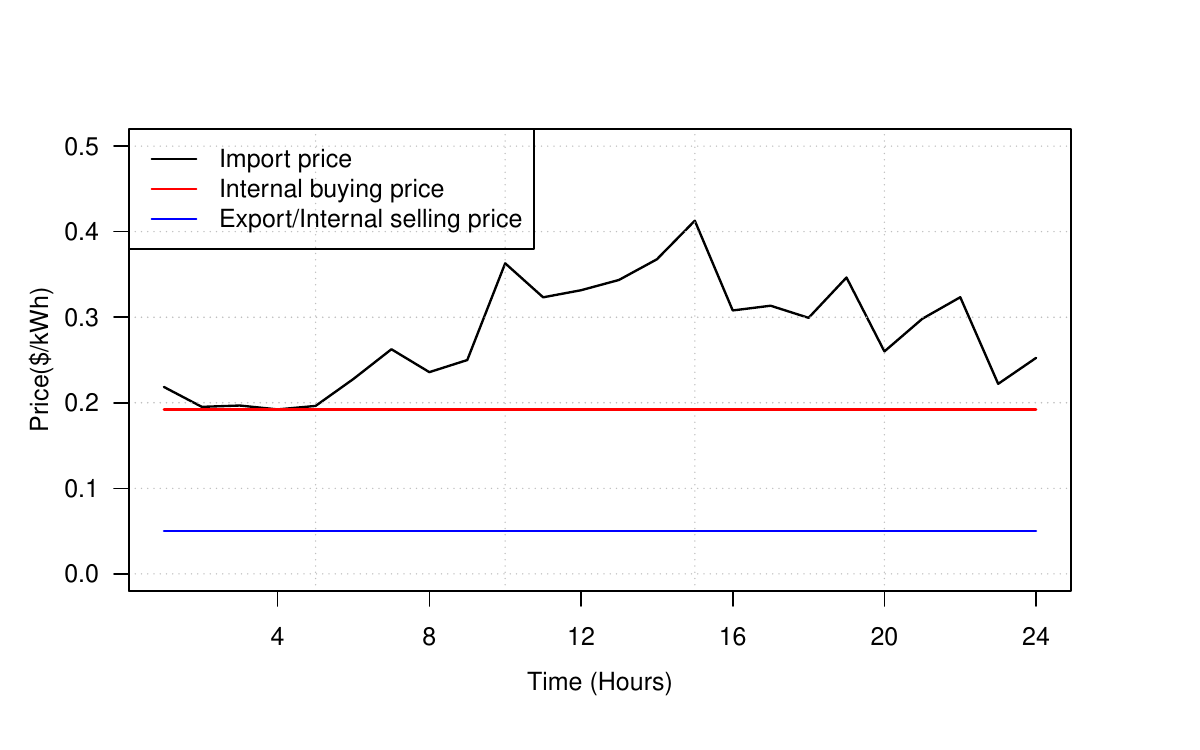}
\caption{The electricity prices for June 18th, 2019 in Cambridge, MA}
\label{fig:price}
\end{center}
\end{figure}

In all three models, we use a time step $\Delta t=1$ hour. The execution period is performed day-ahead from 00:00 till 24:00, where the cycle of the day is divided into $T$ = 24 time slots. We set the total planning period as 10 years and we compute the NPV with a discount rate $r=10\%$.

\subsection{CES Deployment Analysis}
The analysis is done by comparing the strategic decisions corresponding to the different clusters formed across the three models as well as the battery capacity and location decisions and finally the associated costs. Note that for the optimization models applied to this instance, the island model had 4.6\% gap and the interconnected model had a 0.13\% gap. The ESCO model was solved to optimality.

\subsubsection{Batteries Capacity and Location}\label{sec:batterycap}
The capacities and locations of CESs in three %clustering 
\red{business models} are compared in Figure \ref{fig:Batteries}, and the average operations of each CES are also shown in the figure. The numbers on top of the figure are the identifiers of 16 potential battery deployment sites. \red{we assume that the batteries are located at prosumer sites and they are located at 20\% of the prosumers. In this case, as we have 80 prosumers, 16 of these sites are potential battery locations. These 16 locations are selected such that they cover the geographical area uniformly for a given instance.} %\rred
%{\sout{Recall that these 16 potential battery locations are evenly assigned to 80 prosumer's locations, so the displayed numbers represent 16 prosumers' location identifiers.}} 
To examine the efficiency of CESs in three %clustering 
\red{proposed settings,} the average capacities and the maximum capacities are given in Figure \ref{fig:Batteries}. 
Because of the high fuel price in the island model, the demand for energy storage capacity in the island model is relatively high compared to other ones. The number of chosen CESs 
%(but also and their capacities, see Section \ref{sec:batterycap}) 
reached the maximum as we have all 16 batteries assigned for the island model. 
%As explained before, 
On the contrary, for the interconnected model and for the ESCO model, only seven, respectively three, CES are employed. Moreover, the chosen CES capacity for the island models is also bigger. In the interconnected model, the seven CESs are of different capacities based on the number of the prosumers in the cluster. The ESCO \red{model} 
%clustering 
is the most economical setting in terms of the minimum CES capacities required. The results show that the neighborhood requires only three CESs with minimum capacity to meet the community electricity demand. 
Figure \ref{fig:clusters} illustrates how the prosumers and consumers are clustered with assigned CESs to each community for the ESCO model. For this small demonstrated neighborhood, there are 80 prosumers and 40 consumers which are marked in green and blue dots respectively. The red circles are the chosen CES locations out of 16 potential deployment sites (3 batteries in this case) and different colors represent different clusters on the map.

We observe that in all three %clustering 
\red{business models} considered in our study, location B9 and B16 are deployed with batteries which shows that they are needed in these two 
areas.%\footnote{\red{But it could also be that the households from this area are assigned to batteries from the surrounding clusters, correct? Because I think for a household we can choose the assignment to a neighboring battery, correct? }\bissan{Yes but in Figure 6 I see for the IL model we have all CES even 22 and 30 so I am not sure of the results now and I need to go back to the MIP solution}}

\begin{figure}[h]
\includegraphics[trim={0 0 0 2cm},clip,width=1.1\textwidth]{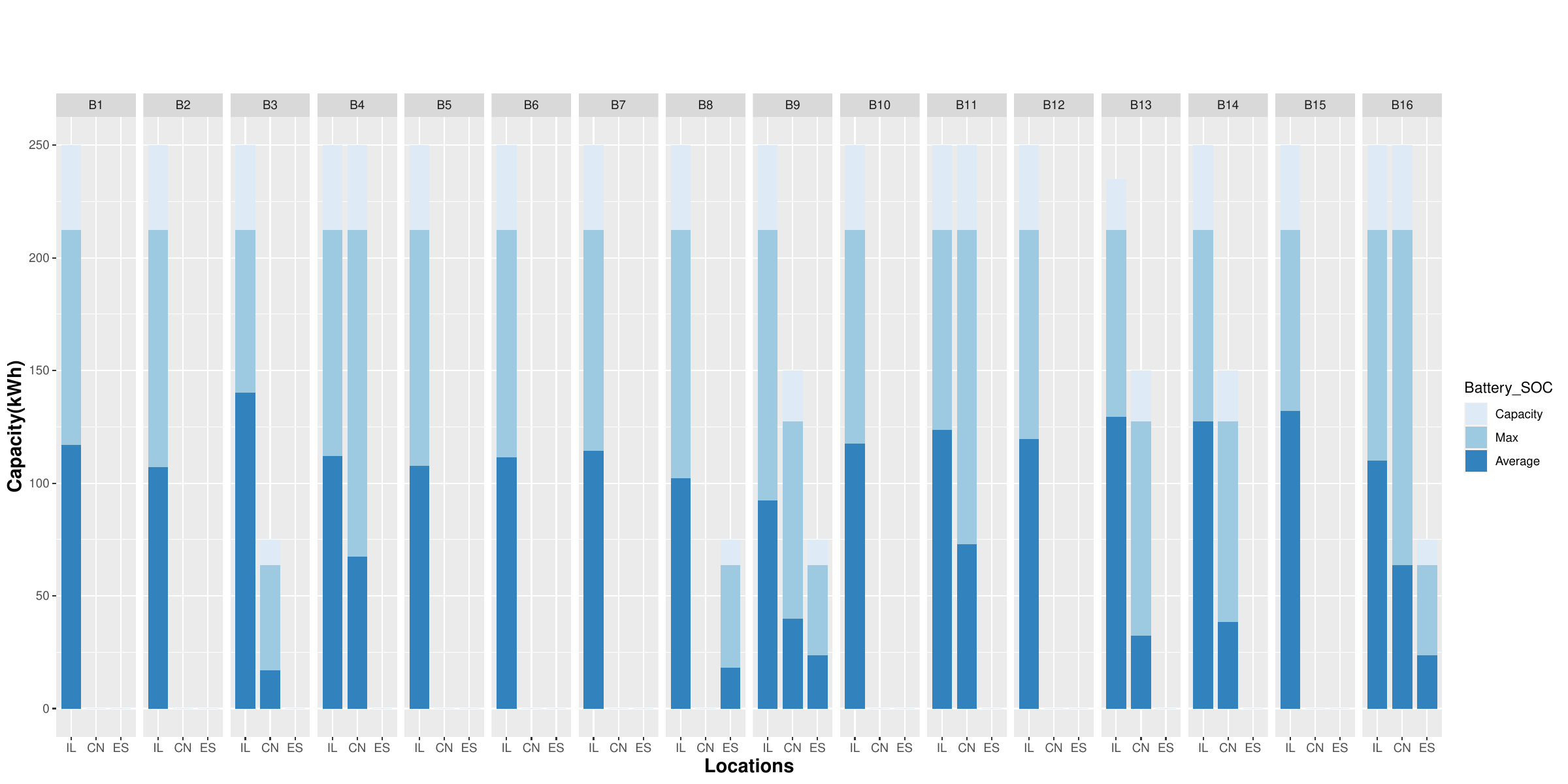}
\caption{Battery deployment at candidate locations  of three clustering methods}
\label{fig:Batteries}
\end{figure}

%\subsubsection{Clusters' Comparison}
%In this part, we analyze the clustering results from the three different models. 
%The interconnected clusters are better scattered and 7 batteries are chosen. For the ESCO model, 3 batteries are assigned and there is no overlap in the clusters. 
%By analyzing the topology of the three solutions, we can derive some additional observations. 
\begin{figure}[h]
%\begin{minipage}{0.3\textwidth}
%\includegraphics[width=1\textwidth]{figs%/I.png}\\
%\centering \small (a) Island model %(\red{16} clusters)
%\end{minipage}
%\hfill
%   \begin{minipage}{0.35\textwidth}
%\includegraphics[width=0.865\textwidth]{figs/C.png}
%\\ \small(b) Interconnected model (6 clusters)
 %  \end{minipage}
 %  \hfill
      \begin{minipage}{\textwidth}\centering
\includegraphics[width=0.5\textwidth, angle =90]{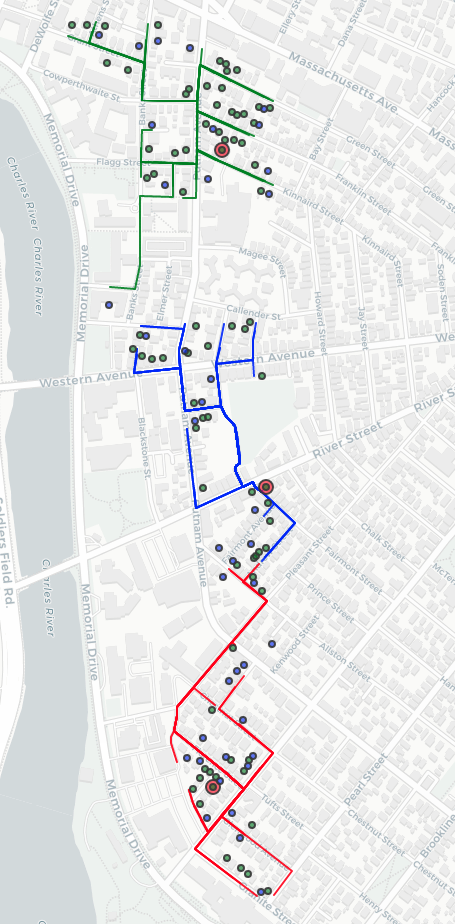}\\
%\centering \small (c) ESCO model (3 clusters)
   \end{minipage}
   \caption{Illustration of CES Locations and Community Networks of the ESCO Model}
   \label{fig:clusters}
\end{figure}

\subsubsection{Batteries Economical Analysis}

The optimal CES configurations of three types of clusters and their corresponding properties in terms of energy sharing efficiency are shown \red{in Table \ref{tab:bat1}.} The district capacity represents the total installed storage capacity in the neighborhood. The number of CESs deployed represents the number of clusters/communities in the neighborhood. The average capacity of CES shows the total energy storage demand of one cluster/community and the household average capacity \red{(i.e., the district capacity divided by the total number of households)} shows the individual energy storage shares when given different clustering objectives.
\begin{table}[H]
\centering
\begin{tabular}{lllll}
\toprule
\multicolumn{2}{l}{Clustering Model}  & Island    & Interconnected & ESCO   \\
\midrule
\multicolumn{2}{l}{Number of CESs deployed}   & 16   & 7     & 3      \\
\multicolumn{2}{l}{District capacity(kWh)} & 4,000  & 1,275  & 225    \\
\multicolumn{2}{l}{Average Capacity(kWh)}  & 250  & 182  & 75     \\
\multicolumn{2}{l}{Household Average Capacity(kWh)}   & 33.3 & 10.6   & 1.9 \\
\multicolumn{2}{l}{Total Battery Dis/Charging(kWh) }  &  2,407.7 & 956.3 & 168.8\\
\multicolumn{2}{l}{Average Battery Dis/Charging(kWh) }  &  171.9  & 136.6 & 56.2\\
\multicolumn{2}{l}{Average Dis/Charging per Household(kWh)}   & 20.1 & 7.9 & 1.4\\
\bottomrule
\end{tabular}
\caption{Overview of CES configurations in the neighborhood}
\label{tab:bat1}
\end{table}
\iffalse
%\blue{Comment: For the average household capacity, the ESCO is 1.875 which is very low. I do not think this is reasonable. For the Island it is more than double the individual household battery with capacity of 13.5, does each household utilize the full 33.33?. Can we validate these. Are you removing the households that are not connected to any CES from this average? Can we comment on the battery utiization in this case? So maybe average as well as min and max values?\ \\}
%\bissan{we need to comment on the business model of the ESCO as it has a negative profit due to the investment cost, additionally we are forcing at least one CES to be present and that is why the number of CES as well as the Average Dis/Charging per Household(kWh) are very low compared to the other models, the number of households of 120 is not profitable for ESCO to be able to have a ROI of the CES}
\fi
The household average capacity of island model is much bigger than the average capacity in the other two models, as well as than the average single household battery's capacity (which is currently 13.5kWh, see e.g.,  \cite{van2018techno}). This result shows %\red{the inefficiency of the island clustering model. Because of the inflexible energy operation and deficiency of distributed energy generated in the island mode microgrid}\footnote{I am not sure about this sentence}, 
that the total cost including CES deployment and fuel cost is much higher in the island %clustering 
\red{business model.} In ESCO model, because of the main decision maker is the energy service company whose main goal in this energy transaction activity is to maximize its profits, they make good use of the price difference of utility company from the main grid but do not maximize the battery utilisation rate, which is the reason why the average capacity of CES in ESCO model is very low and thus the energy storage utilisation is not efficient. In terms of economical effect from the CES clustering, Table \ref{tab:bat2} compares the electricity bill of the neighborhood after applying three different %clustering 
models:
\begin{table}[H]
\centering
\begin{tabular}{lllll}
\toprule
\multicolumn{2}{l}{%Clustering 
\red{Business} Model}                               & Island    & Interconnected & ESCO   \\
\midrule
\multicolumn{2}{l}{Total CES investment(k\$)}                                              & 1,720        & 551    & 99 \\
\multicolumn{2}{l}{Total MT investment(k\$)}                                              & 250        & N/A    & N/A \\
\multirow{2}{*}{\begin{tabular}[c]{@{}l@{}}Total Operational\\ Cost (k\$)\end{tabular}}     
                                                                        & Fuel Cost        & 4,337        & N/A       & N/A     \\
                                                                  %      & Internal Cost    & N/A         & N/A       & 246 \\
                                                                        & Revenue & N/A            & 200      &616  \\
                               & Cost    & N/A          & 3,779     & 177
 \\
                                                               %     & External Revenue & N/A          & 200       & 304  \\
%                                                                      & Net Cost         & 4,818        & 3,179     & 4,292 \\ 
\multicolumn{2}{l}{Total cost(k\$)}                                                        & 6,307        & 4,130     & -340 \\  
\bottomrule
\end{tabular}
\caption{Overview of households electricity cost in one neighborhood of 120 households}
\label{tab:bat2}
\end{table}
\iffalse
{need to explain what is meant by Internal Cost and Revenues in the fixed cost section for the ESCO? report the gaps \ \\}\fi
In a planning period of 10 years, the total cost of the energy sharing neighborhood consists of two parts: the initial CES/MT deployment investment and the total operational cost. In the island model, the operational cost is only from the fuel consumption of microturbines while in the interconnected setting, the cost consists of the cost of importing electricity from the main grid minus the profit for the prosumers who sell the energy they produce. In the ESCO setting, there is the internal transaction between prosumers and consumers which is operated by ESCOs. The internal cost is the electricity bill paid by households who get electricity from CES and the internal revenue is the profit that prosumers get from charging the CESs. As we can see in the table, the {total energy cost over the planning period of 10 years is the highest for the island model, and it is 1.5 higher than the cost of the interconnected model.} %Even though} the ESCO model is the most economical in terms of total cost, this operation mechanism may pose a significant challenge to the operability of main grid  because of the highly intermittent solar energy generation of prosumers. \red{For the ESCO model, although the model maximizes the profits however when taking the investment costs into account they are higher compared to the operational revenue.}\footnote{Also, in the table total CES investment says s 496, in the text it says 401.

\subsection{Daily Operation Analysis}
In addition to the strategic long-term differences in the three models, in this section we present the operational differences in the households' as well as batteries' operations.
\subsubsection{Batteries' Daily Operation}
The daily electricity charging operation of one arbitrary CES in three different %clustering 
\red{business models} is shown in Figure \ref{fig:SOC}. As we can see in Figure \ref{fig:SOC}, the CESs in the island and the interconnected model follow the same charging and discharging pattern, which are positively correlated to the amount of energy that prosumer produce (radiation rate). However, in the ESCO model, the charging and discharging operation is manipulated by the energy service company to maximize their profits. %\footnote{\bissan{comment on investment costs}}%, \red{that is to say they discharge the CES once they get enough energy then sell it to the households in the community. In term of the battery utilization efficiency, the CESs in the ESCO model are not economically efficient.}\footnote{I dont understand it}
\begin{figure}[H]\begin{center}
\includegraphics[width=0.8\textwidth]{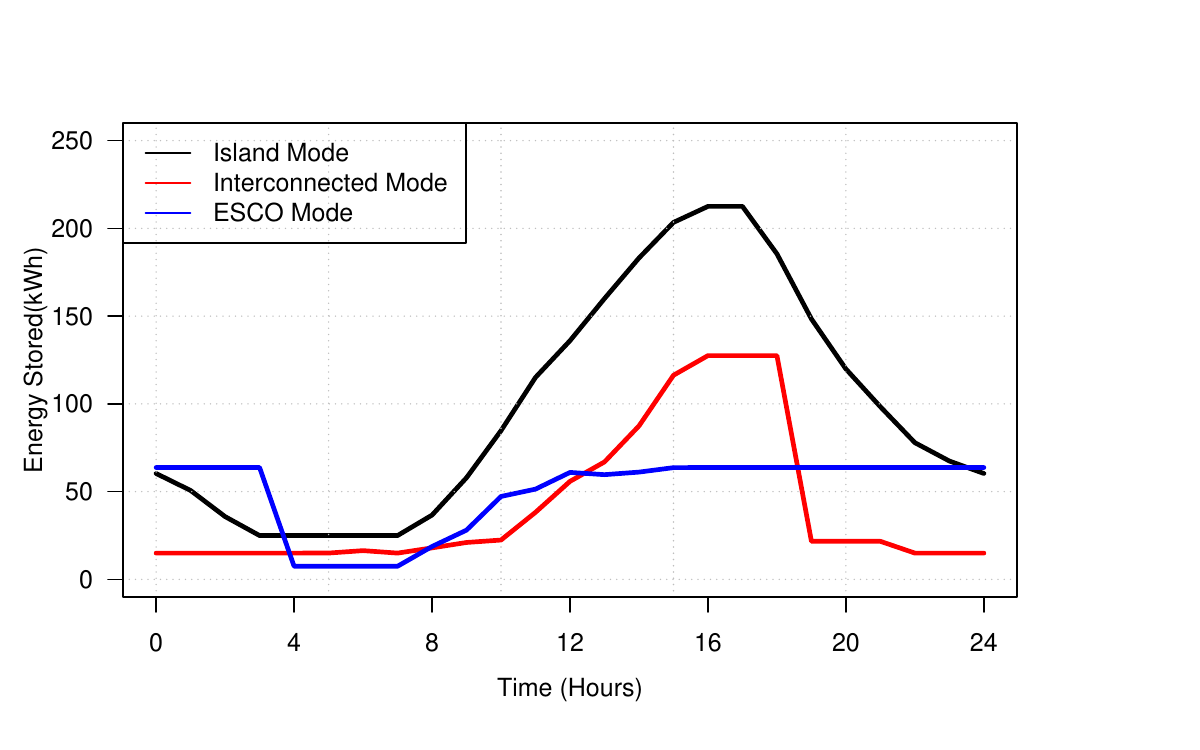}
\caption{Daily SOC operations of each battery of three %clustering methods
\red{considered business models}}
\label{fig:SOC}\end{center}
\end{figure}
\subsubsection{Households' Daily Operation}
Figures \ref{fig:flow1} - \ref{fig:flow3} show the daily household energy utilization curve of one arbitrary consumer and one arbitrary prosumer for the corresponding three models. Straight lines represent the external energy exchange with the main grid or microturbine while the dotted lines represent the internal energy exchange between consumers and prosumers through CESs. Black lines are referred to the consumers while the red lines are referred to the prosumers. In the island model, the internal energy exchange operation \red{is higher than that of the other models}, because of the high fuel price used by microturbine. The CES charging curve in the island model is in accordance with the energy generation of prosumers. In the other two %clustering 
\red{business models,} this feature is not obvious because of external energy exchange with main grid. The energy exportation peak time to the main grid also varies in interconnected and ESCO model. The energy consumption for both consumers and prosumers is mainly provided by main grid at the peak period, whereas in the island model, the electricity need is satisfied primarily by the CES. {When the electricity stored in the CES cannot satisfy the households' load, they start using electricity from the microturbines.}%\footnote{I dont understand this}

\begin{figure}[H]
\begin{minipage}{0.45\textwidth}
\includegraphics[width=1.18\textwidth]{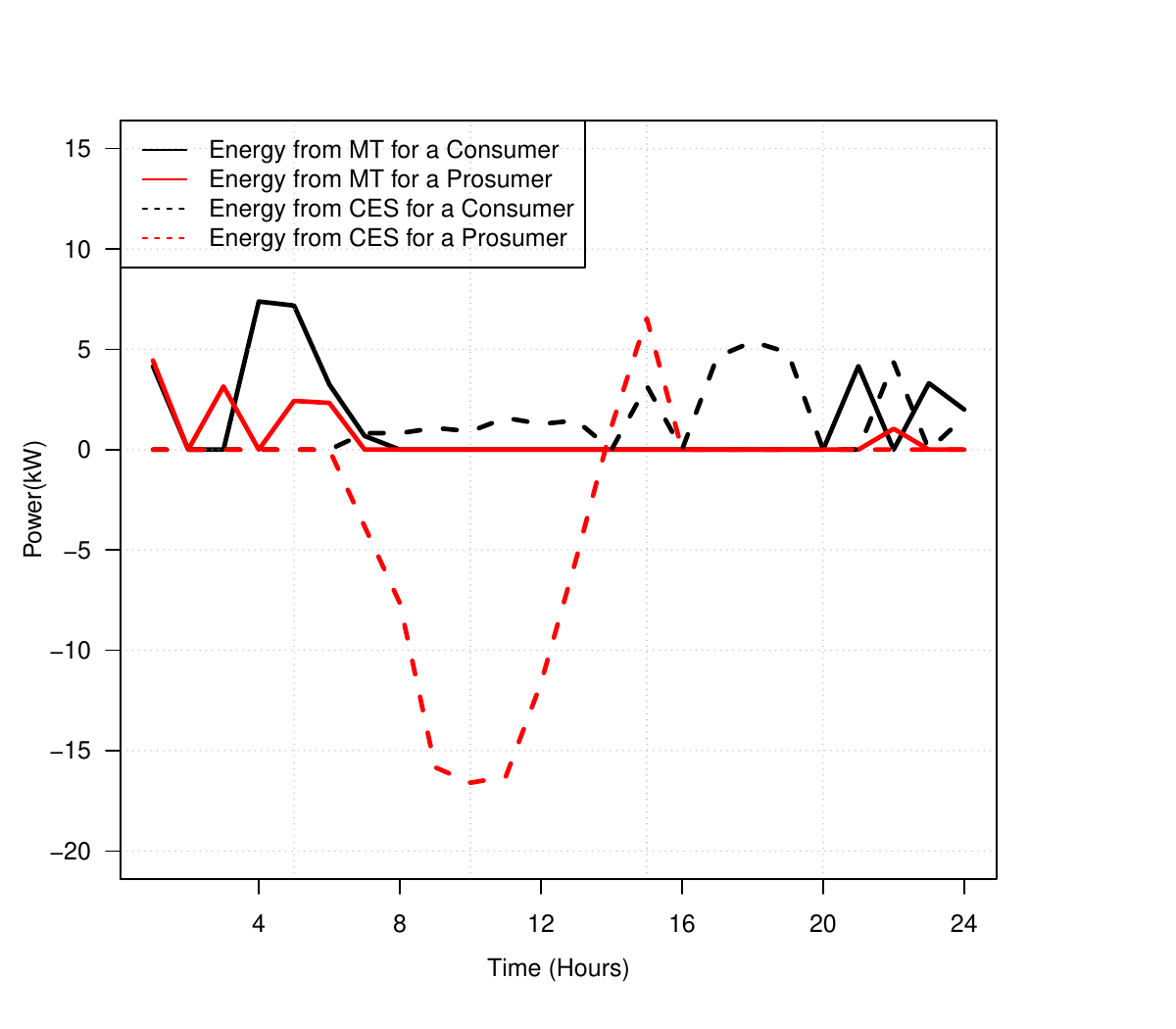}
\caption{An arbitrary consumer and prosumer's daily power flow operations in island model}
\label{fig:flow1}
\end{minipage}\qquad
\begin{minipage}{0.45\textwidth}
\includegraphics[width=1.18\textwidth]{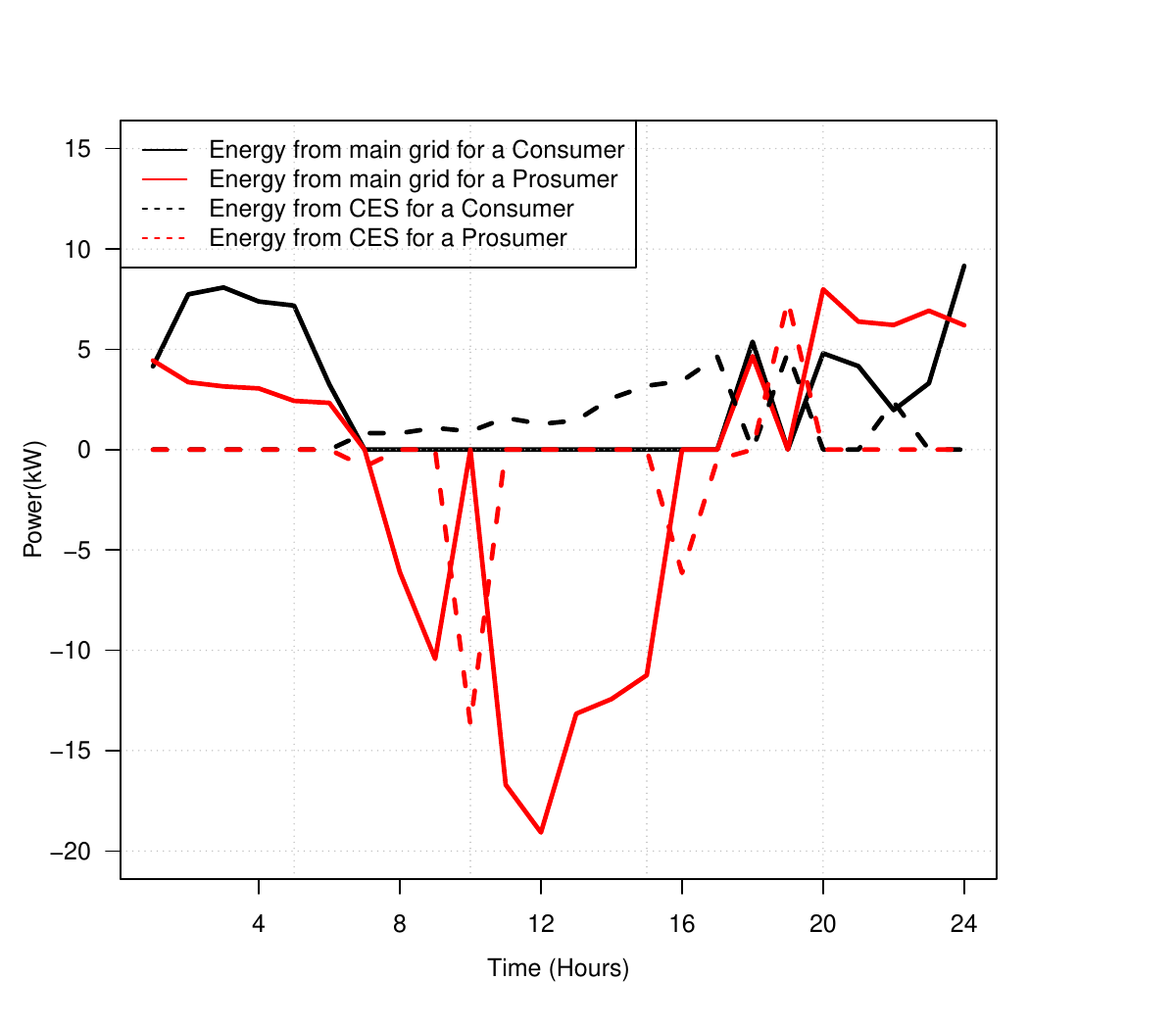}
\caption{An arbitrary consumer and prosumer's daily power flow operations in interconnected model}
\label{fig:flow2}
\end{minipage}
\begin{center}
\begin{minipage}{0.45\textwidth}
\includegraphics[width=1.18\textwidth]{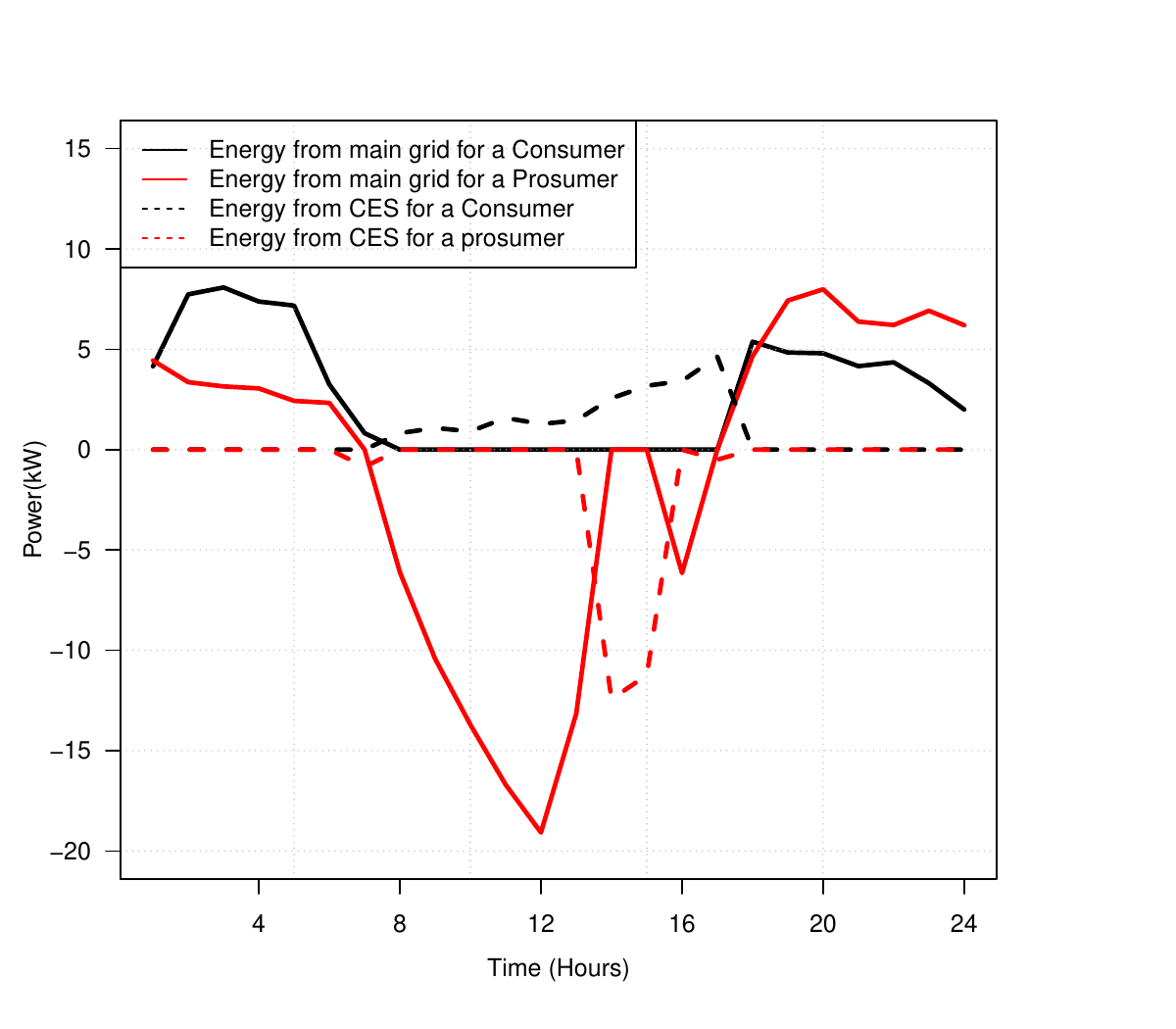}
\caption{An arbitrary consumer and prosumer's daily power flow operations in ESCO model}
\label{fig:flow3}
\end{minipage}\end{center}
\end{figure}

\subsection{Clustering Benefits}

To compare the advantage of CES to HES, we also simulated the daily optimal operations for an arbitrary consumer and arbitrary prosumer when equipped with a HES at a prosumer's house. Figure \ref{fig:HESflow} shows the optimal daily power flow operations. We assume that the Household Energy Storage installed at each prosumer's house %with the optimal 
{has a capacity of HES 13.5 kWh, as assumed in} \citep{barbour2018community}. Because of the power limit from the HES setting, the excessive energy generated from the solar panel goes directly to the main grid which makes the self-consumption not as good as when prosumers share CESs.

\begin{figure}[H]\begin{center}
\includegraphics[width=0.7\textwidth]{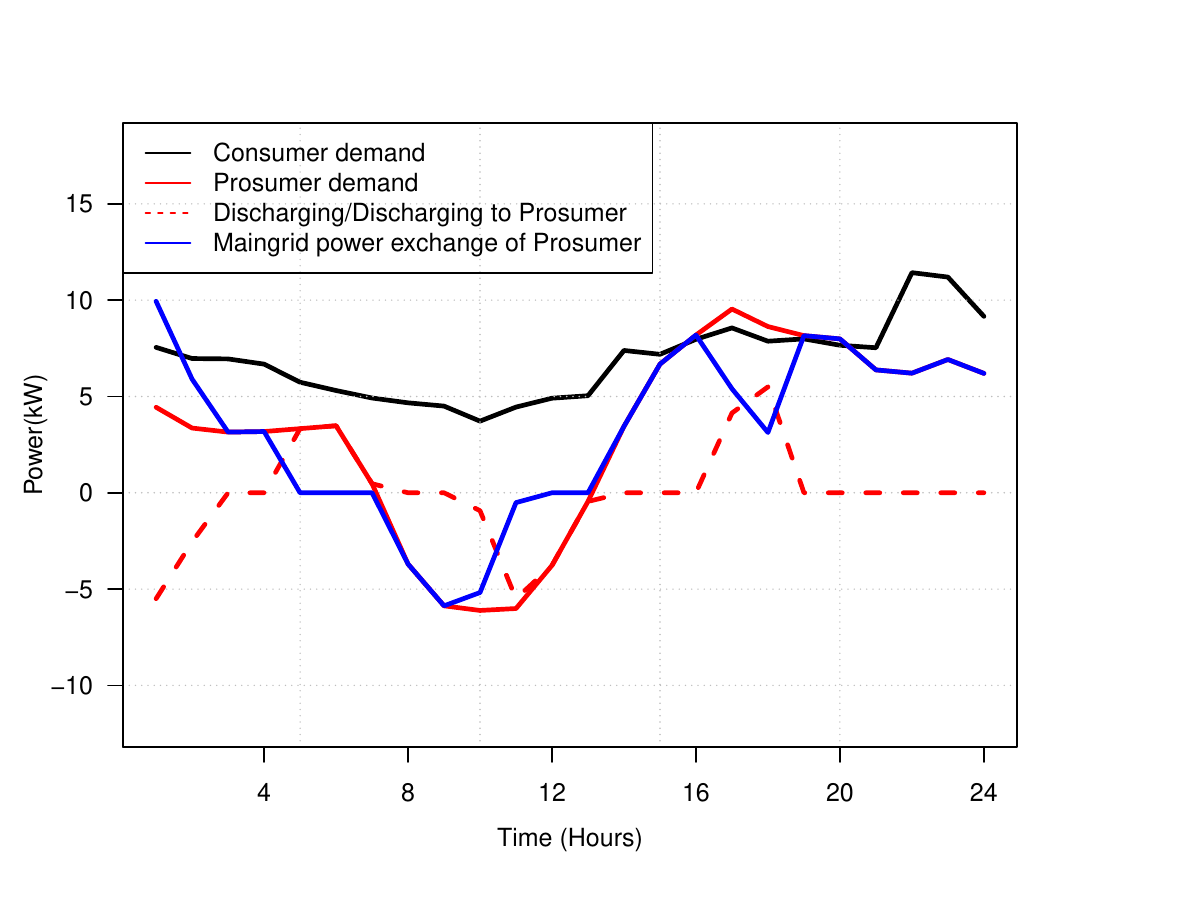}
\caption{An arbitrary consumer and prosumer's daily power flow operations with HES}
\label{fig:HESflow}    
\end{center}
\end{figure}
The price of single HES installed at each prosumer's house is 15,175\$ from Sonnen company \cite{SonnenPrice}, thus the total battery cost {for deploying 80 HES at prosumer locations} 
%of this neighborhood 
is 1.214 M\$ (there is no electricity generation at consumers' houses so they are not equipped with HES), which is much higher than the CES deployment cost in interconnected and ESCO clustering neighborhood (551 k\$ and 99 k\$, respectively, cf.\ Table \ref{tab:bat2}). In terms of electricity bill saving, the total electricity bill in this neighborhood without any energy storage is 4,338 k\$. In comparison, for the interconnected model total cost for households is around 3,579 k\$, which shows that community energy storage can decrease the community electricity bill significantly. On the other hand, because of high fuel consumption cost and big capacity CES deployed in the island model, this clustering model cannot be considered economical from the current perspective.%\footnote{We do not comment ESCO model here}
{
\subsection{Sensitivity Analysis}
In this section, we examine the sensitivity of the results to
the variation in the households' demand as well the variation of production data for prosumer households. To do so, we consider as the “nominal scenario” the one that corresponds to the given consumer and prosumer profiles that we used in 
our case study in Section \ref{sec:data_gen}. For the nominal solution (location and capacity of shared batteries
and assignment of households) obtained from our models, we then reevaluate the operational
expenses for an additional set of 365 different scenarios in which we sample the demand and
production values. These 365 scenarios represent different days throughout the year. We sample these values from a uniform distribution within the interval of $\pm 10\%$ of the nominal value. We present the variation in daily operational costs in Figure \ref{fig:sensitivity} with the operational costs of the “nominal scenario” shown in the dotted line. For the ESCO model, the operational profit is given.

The obtained results show that the operational costs/profits for the three different models vary within a very small range. In particular, for the island model the variation is around 5\% while for the Interconnected and the ESCO models the variation is 1\% and 2\%, respectively. Therefore, our study indicates that the proposed models are robust towards changes in the demand.
Indeed, the operational costs are not significantly impacted by the fluctuation in the demand happening throughout the year.  Additionally, the total operational costs of the Island model is less than 10\% of the total costs (9.75\%) and the variation of the daily operational costs (over the planning period) is 0.086\% of the total costs. For the Interconnected model, the total operational costs form 13.6\% of the total costs and the variation of the daily operational costs is 0.057\% of the total costs. Finally, for the ESCO model, the total operational profit is 18.5\% of the total profit and the variation of the daily operational profit is 0.051\% of the total profit.
\begin{figure}[H]
\includegraphics[width=1\textwidth]{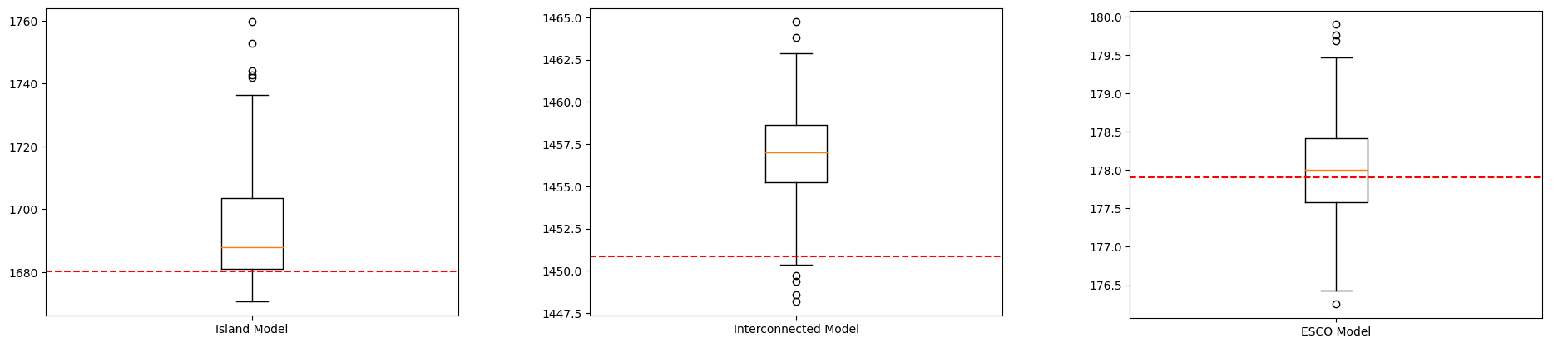}
\caption{Daily operational cost (profit for the case of ESCO) variation using 365 demand scenarios for the three models.}
\label{fig:sensitivity}

\end{figure}

}
\section{Computational Performance}\label{sec:results}
{In this section we studied the computational performance of the proposed mathematical models, with the aim of exploring \red{their limitations,} in terms of size of the input data for which (nearly)-optimal solutions can be found.} The goal is to test the effectiveness and scalability  of the proposed \red{formulations} and the associated Benders decomposition approach.
Our computational experiments have been performed on a workstation computer with an Intel i5-1035G1 1.0 GHz CPU with 12 GB of RAM, running Windows 10 version 21H2. IBM CPLEX %{\textregistered} 
20.1.0.0 (referred to as CPLEX in the following) was used as a general-purpose MILP solver.  For all reported runs we set the time limit to 3600s. CPLEX parameters are set to their default values (e.g.,  the relative MILP gap tolerance is 0.01\%).
Benders decomposition approach is implemented using the annotated Benders strategy in CPLEX, which allows to efficiently exploit some of the state-of-the art computational techniques including cut-loop separation, stabilization using in-out strategy, multi-start and internal handling of data structures (see e.g., \cite{cplex-benders} for more details).  

\paragraph{Benchmark Instances} For this experiment we have generated two types of benchmark instances based on the ratio between the prosumers and the consumers in the network:  in the first group we set $|P|=2|C|$, and in the second group we set $|P|=3|C|$. The number of consumers $|C|$ is chosen from the set $\{30,40,50,75,100,150,200,250,300\}$, and the number of prosumers is set correspondingly. Detailed computational results are provided in Tables \ref{tab:I1}-\ref{tab:G2} where we report: 
the number of consumers $|C|$,  the number of prosumers $|P|$,
the type  of algorithm we used (where  
com  refers to the compact model,
and   Ben  to the Benders decomposition approach),
 the number of nodes processed in the  branching tree (\#Nodes),
 %the value of the LP-relaxation at the root node ($LB_{root}$), 
 the value of the global dual bound at the end of the run ($DB$), 
the gap (in percentage) between the global dual bound and the best integer solution {found by the corresponding approach} ($Gap[\%]$),
the value of the optimal (or best-found) integer solution ($\ObjVal$),
the CPU time in seconds (CPU$[s]$).  TL means that the time limit was reached for this instance) and if CPLEX was unable to solve the LP model, we denote it by $-$ in the table. \red{The default time limit was set to two hours.}

\subsection{Computational Performance for the Island Model}

Table \ref{tab:I1} reports the obtained results for instances with up to 300 consumers and 600 prosumers. We observe that both computational approaches have difficulties in finding optimal solutions for the island  model. The final gaps are below 20\% for instances with up to 600 households ($|P| + |C|$) in total. {Instances above 600 households could not be solved by the compact model. Recall that the associated compact MILP formulations contain around millions of variables and constraints. Hence, it is not surprising that 
%the problems are difficult when it comes to proving the optimality. 
this poses serious limitations to the MILP solver, so that even solving the LP-relaxation cannot be done within two hours. In these situations, Benders decomposition provides significant advanatges over solving the compact model - even for the largest instances from this data set with up to 900 households, Benders decomposition manages to find feasible dual and primal bounds, with reasonable gaps that remain below 15\%.}    

When it comes to computational performance for instances with up to 200 consumers and 400 prosumers, there is no clear dominance between the two proposed approaches. Better primal and dual bounds are found by different methods across the instances. We observe that Benders decomposition draws the advantage of having a smaller number of variables, so that the branching process starts earlier and many more branching nodes could be explored within the time limit. On the other hand, this advantage does not necessarily translate into a better quality of lower bounds, which can be explained by the fact that CPLEX can strengthen the lower bounds by deriving general-purpose constraints from the full model description, which is provided by the compact formulation. However, the main advantage of using Benders in this case is that one can solve larger instances (up to 900 households) with a reasonable gap while the compact model is not able to solve these instances.
 
For Table \ref{tab:I2}, instances up to 600 households can be solved by both approaches within the time limit. Solving larger instances 
%suffered from 
was prohibitive due to memory limitations.

\begin{table}[h!]
\centering
%\scriptsize
\tabcolsep 7pt
\renewcommand
\arraystretch{1.2}
\begin{tabular}{lrrrrrr}
\toprule			
$|C|$ & $|P|$ &  Alg  %&  \#Cuts  
&  \#Nodes 
%& $LB_{root}$ 
& DB  & Gap$[\%]$ & $\ObjVal$ \\ % &  CPU$[s]$ \\
	\cline{1-2} 
	\cline{4-7} %\cline{4-8}
	%\cline{6-9}
		
30	& 60	& com	&  1192	& 4,226,289.59	& 2.78	& 4,347,052.52\\ % 	& TL \\
40	& 80	& com	& 10	& 5,826,504.22	& 6.87	& 6,256,339.30\\ % 	& TL  \\
50	& 100	& com	& 0	& 7,503,323.47	& 6.32	& 8,009355.67\\ % 	& TL\\
75	& 150	& com	& 0	& 11,044,866.51	& 17.81	& 13,437,821.48\\ % 	& TL\\
100	& 200	& com	& 0	& 13,601,169.84	& 14.39	& 15,887,738.74\\ % 	& TL\\
150	& 300	& com	& 0	& 20,096,716.59	& 14.26	& 23,438,927.89\\ % 	& TL\\
200	& 400	& com	& 0	& 27,268,151.08	& 13.34	& 31,464,721.68\\ % 	& TL\\
250	& 500	& com	& -	& -	& -	& 40,739,926.39\\ % 	& TL\\
300	& 600	& com	& -	& -	& -	& 47,780,669.55\\ % 	& TL\\
\hline
30	& 60	& Ben	& 40815	& 4,259,386.63	& 8.81	& 4,670,811.89\\ % 	& TL \\
40	& 80	& Ben	& 11969	& 5,821,760.12	& 17.30	& 7,039,600.17\\ % 	& TL\\
50	& 100	& Ben	& 2987	& 7,483,812.28	& 18.40	& 9,171,189.75\\ % 	& TL\\
75	& 150	& Ben	& 326	& 11,017,784.36	& 18.01	& 13,437,821.48\\ % 	& TL\\
100	& 200	& Ben	& 0	& 13,568,964.44	& 14.59	& 15,887,738.74\\ % 	& TL\\
150	& 300	& Ben	& 0	& 20,049,290.47	& 14.46	& 23,438,927.89\\ % 	& TL\\
200	& 400	& Ben	& 0	& 27,268,144.31	& 13.34	& 31,464,721.68\\ % 	& TL\\
250	& 500	& Ben	& 0	& 34,958,988.54	& 14.19	& 40,739,926.39\\ % 	& TL\\
300	& 600	& Ben	& 0	& 40,785,761.00	& 14.64	& 47,780,669.55\\ % 	& TL\\
\bottomrule
\end{tabular}
\caption{Results for the Island Model with $|P|=2|C|$ 
%and time limit of \red{two hours}
%one hour
\label{tab:I1}}
\end{table}

\begin{table}[h!]
\centering
%\scriptsize
\tabcolsep 7pt
\renewcommand
\arraystretch{1.2}
\begin{tabular}{lrrrrrr}
\toprule			
$|C|$ & $|P|$ &  Alg  %&  \#Cuts  
&  \#Nodes 
%& $LB_{root}$ 
& DB & Gap$[\%]$ & $\ObjVal$ \\ % &  CPU$[s]$ \\
	\cline{1-2} 
	\cline{4-7}
30	&90	&com &2236&	3,467,549.65&	0.13&	3,472,182.77  \\ %&	TL \\
40	&120	&com	&0	&7,880,104.72	&14.32	&9,197,362.93 \\ %	&TL\\
50	&150	&com	&0	&9,882,254.79	&15.54	&11,700,811.7 \\ %	&TL\\
75	&225	&com	&0	&13,285,520	&16.56	&15,922,081.87 \\ %	&TL\\
100	&300	&com	&0	&17,319,043.76	&15.92	&20,597,863.94 \\ %	&TL\\
150	&450	&com	&0	&26,617,570.61	&14.94	&31,294,052.47 \\ %	&TL\\

\hline
30	&90	&Ben	&7757	&5,682,616.73	&16.56	&6,810,316.41 \\ %	&TL\\
40	&120	&Ben	&808	&7,897,619.18	&17.17	&9,535,106.10 \\ %	&TL \\
50	&150	&Ben	&285	&9,835,027.06	&17.96	&11,987,958.32 \\ %	&TL\\
75	&225	&Ben	&0	&13,285,540.43	&16.56	&15,922,081.98 \\ %	&TL\\
100	&300	&Ben	&0	&17,319,059.19	&15.92	&20,597,863.94 \\ %	&TL\\
150	&450	&Ben	&0	&26,558,702.20	&15.13	&31,294,052.47 \\ %	&TL\\
\bottomrule
\end{tabular}
\caption{Results for the Island Model with $|P|=3|C|$
%and time limit of \red{two hours}
%one hour
\label{tab:I2}}
\end{table}

\subsection{Computational Performance for the Interconnected Model}
We now turn our attention to the interconnected model. We consider instances with up to 300 consumers and 900 prosumers. Detailed results are provided in Tables \ref{tab:C1}-\ref{tab:C2}.
We observe that for the instances of the same size, the computational performance of the compact model and the Benders reformulation is much better when compared to the Island Model. Most of the instances with up to 150 consumers are solved to optimality (or within a small optimality gap, particularly using the Benders approach). \red{In Table \ref{tab:C1} we do not report the computing times, since the time limit was reached for all except two instances, for which the numbers given in brackets in the gap column correspond to the computing times needed to achieve this gap. }

This time we notice a significant difference between compact model and Benders reformulation which is particularly striking when considering larger instances. For the setting with $|P|=2|C|$ the compact formulation reaches its limits with 200  consumers -- the gap of 24\% is achieved for $|C|=200$ and the LP-relaxation could not be solved within two hours for the larger instances. The same effect is even more pronounced for instances with   
$|P|=3|C|$ where the compact formulation reaches its limits already with 150 consumers and the gap is much larger compared to Benders (22.85\% versus 1.61\%). Indeed, for none of the instances with 200 or more consumers, the LP-relaxation could be solved within two hours. On the contrary, Benders reformulation performs quite stable and scales well with the size of the instances. While the final gaps increase with the size of the input data (due to the fact that less branching nodes can be enumerated within the time limit), the quality of primal and dual bounds remains satisfactory with the largest gaps remaining below 7\% for 1000 households and reaching 24\% for 1200 households.  

\begin{table}[h!]
\centering
%\scriptsize
\tabcolsep 4pt
\renewcommand
\arraystretch{1.2}
\begin{tabular}{lrrrrrr}
\toprule			
$|C|$ & $|P|$ &  Alg  %&  \#Cuts  
&  \#Nodes 
%& $LB_{root}$ 
& DB  & Gap$[\%]$ & $\ObjVal$  \\ % &  CPU$[s]$ \\
	\cline{1-2} 
	\cline{4-7}
	%\cline{6-9}
30	& 60	& com	& 20200	& 2,617,407.51	& $\leftidx{^{(6493s)}}{0.01}{}$  	& 2,617,656.25 \\ %	& 6493.35 \\
40	& 80	& com	& 93	& 3,574,282.48	& 0.14	& 3,579,307.24 \\ %	& TL\\
50	& 100	& com	& 1334	& 4,627,122.19	& 0.27	& 4,639,547.08 \\ %	& TL\\
75	& 150	& com	& 0	& 6,779,058.61	& 0.07	& 6,784,090.61 \\ %	& TL\\
100	& 200	& com	& 0	& 8,266,172.87	& 3.25	& 8,543,549.08 \\ %	& TL\\
150	& 300	& com	& 0	& 12,328,748.14	& 24.58	& 16,347,530.57 \\ %	& TL\\
200	& 400	& com	& 0	& 16,767,702.86	& 24.06	& 22,080,661.18 \\ %	& TL\\
250	& 500	& com	& 0	& -	& -	& 28,351,946.06 \\ %	& TL\\
300	& 600	& com	& 0	& -	& -	& 33,102,698.18 \\ %	& TL\\

\hline
30	& 60	& Ben	& 28401	& 2,617,397.46	& 
$\leftidx{^{(633s)}}{0.01}{}$  	& 2,617,656.25 \\ %	& 633.33\\
40	& 80	& Ben	& 9360	& 3,533,774.82	& 1.36	& 3,582,515.98 \\ %	& TL\\
50	& 100	& Ben	& 5849	& 4,550,002.47	& 1.92	& 4,638,955.57 \\ %	& TL\\
75	& 150	& Ben	& 11441	& 6,713,158.55	& 1.05	& 6,784,090.61 \\ %	& TL\\
100	& 200	& Ben	& 1855	& 8,249,879.63	& 0.43	& 8,285,296.57 \\ %	& TL\\
150	& 300	& Ben	& 2299	& 12,243,800.8	& 0.87	& 12,351,782.79 \\ %	& TL\\
200	& 400	& Ben	& 414	& 16,468,563.96	& 2.41	& 16,875,265.24 \\ %	& TL\\
250	& 500	& Ben	& 989	& 21,135,028.23	& 2.29	& 21,630,628.08 \\ %	& TL\\
300	& 600	& Ben	& 0	& 24,542,614.52	& 25.86	& 33,102,698.18 \\ %	& TL\\

\bottomrule
\end{tabular}
\caption{Results for the Interconnected  Model with $|P|=2|C|$ 
%and time limit of \red{two hours}
%one hour
\label{tab:C1}}
\end{table}

\begin{table}[tbh!]
\centering
%\scriptsize
\tabcolsep 7pt
\renewcommand
\arraystretch{1.2}
\begin{tabular}{lrrrrrr}
\toprule			
$|C|$ & $|P|$ &  Alg  %&  \#Cuts  
&  \#Nodes 
%& $LB_{root}$ 
& DB  & Gap$[\%]$ & $\ObjVal$ \\ % &  CPU$[s]$ \\
	\cline{1-2} 
	\cline{4-7}
	%\cline{6-9}
30	&90	&com	&2236	&3,467,549.65	&0.13	&3,472,182.77 \\ %	&TL\\
40	&120	&com	&159	&4,845,690.94	&0.3	&4,860,384.63 \\ %	&TL\\
50	&150	&com	&0	&6,022,186.42	&0.06	&6,025,821.25 \\ %	&TL\\
75	&225	&com	&0	&8,097,850.06	&8.97	&8,896,068.25 \\ %	&TL\\
100	&300	&com	&0	&10,570,748.41	&23.6	&13,836,052.14 \\ %	&TL\\
150	&450	&com	&0	&16,329,693.67	&22.85	&21,165,579.80 \\ %	&TL\\
200	&600	&com	&-	&-	&-	&- \\ %	&TL\\
250	&750	&com	&-	&-	&-	&- \\ %	&TL\\
300	&900	&com	&-	&-	&-	&- \\ %	&TL\\
\hline
30	&90	&Ben	&90247	&3,465,344.34	&0.20	&3,472,182.77 \\ %	&TL\\
40	&120	&Ben	&3572	&4,766,969.12	&1.72	&4,850,292.57 \\ %	&TL\\
50	&150	&Ben	&19374	&5,981,790.53	&0.73	&6,025,821.24 \\ %	&TL\\
75	&225	&Ben	&7243	&7,945,517.55	&1.93	&8,101,939.90 \\ %	&TL\\
100	&300	&Ben	&6860	&10,389,134.68	&1.77	&10,575,873.66 \\ %	&TL\\
150	&450	&Ben	&1572	&16,087,496.27	&1.61	&16,350,398.03 \\ %	&TL\\
200	&600	&Ben	&1020	&22,034,039.70	&3.07	&22,732,957.85 \\ %	&TL\\
250	&750	&Ben	&0	&27,441,200	&6.33	&29,294,900	 \\ %&TL\\
300	&900	&Ben	&0	&33,657,100	&24.01	&44,291,000	 \\ % &TL\\
\bottomrule
\end{tabular}
\caption{Results for the Interconnected  Model with $|P|=3|C|$
%and time limit of \red{two hours}
%one hour
\label{tab:C2}}
\end{table}

\newpage
\subsection{Computational Performance for the ESCO Model}

%\newpage
%As we can observe from the tables above, the unseparable benders algorithm always outperforms the separable benders algorithm especially especially solving big instances, this is because of the fixed continuous variable $S_i^t$ in master problem makes it difficult to find the optimal solutions in subproblem.  
In Tables \ref{tab:G1}-\ref{tab:G2}, we report results of the ESCO model with up to 300 consumers and 600 prosumers. In this case, the compact model also struggles to find high quality solutions for large instances as only instances with up to 100 consumers and 200 prosumers were solved with a reasonable gap while larger instances of 150 consumers or more had 
%an unreasonable 
huge gaps or could not be solved. For the Benders case, instances with up to 300 consumers and 600 prosumers (for $|C|=2|P|$) and up to 200 consumers (for the setup with $|C|=3|P|$) turned to be tractable. The maximum gap reached was less than 25\% for the largest instance from this dataset.
\begin{table}[tbh!]
\centering
%\scriptsize
\tabcolsep 4pt
\renewcommand
\arraystretch{1.2}
\begin{tabular}{lrrrrrrr}
\toprule			
$|C|$ & $|P|$ &  Alg  %&  \#Cuts  
&  \#Nodes 
%& $LB_{root}$ 
& DB  & Gap$[\%]$ & $\ObjVal$ &  CPU$[s]$ \\
	\cline{1-2} 
	\cline{4-8}
	%\cline{6-9}
30	& 60	& com	& 0	& 197,373.12	& 0	& 197,373.12	& 38.76 \\
40	& 80	& com	& 0	& 339,899.10	& 0	& 339,899.10	& 68.49\\
50	& 100	& com	& 0	& 440,332.76	& 0	& 440,332.76	& 61.34\\
75	& 150	& com	& 0	& 667,753.46	& 0	& 667,753.46	& 321.38\\
100	& 200	& com	& 0	& 864,226.55	& 0	& 864,226.55	& 992.72\\
150	& 300	& com	& 0	& 1,242,700.15	& 561.95	& 187,732.09	& TL\\
200	& 400	& com	& 0	& 1,650,325.11	& 2807.03	& 56,770.12	& TL\\
250	& 500	& com	& -	& -	& -	& -	& -\\
300	& 600	& com	& -	& -	& -	& -	& -\\
\hline
30	& 60	& Ben	& 1012	& 197,373.12	& 0	& 197,373.12	& 141.12\\
40	& 80	& Ben	& 3264	& 339,899.10	& 0	& 339,899.10	& 233.91\\
50	& 100	& Ben	& 302	& 440,332.76	& 0	& 440,332.76	& 163.85\\
75	& 150	& Ben	& 22920	& 697,440.88	& 4.45	& 667,753.46	& TL\\
100	& 200	& Ben	& 532	& 864,226.55	& 0	& 864,226.55	& 1118.71\\
150	& 300	& Ben	& 18611	& 1,270,880.68	& 4.32	& 1,218,213.39	& TL\\
200	& 400	& Ben	& 13793	& 1,720,784.13	& 10.54	& 1,556,696.90	& TL\\
250	& 500	& Ben	& 1954	& 2,218,491.63	& 13.63	& 1,952,433.67	& TL\\
300	& 600	& Ben	& 3	& 2,586,454.06	& 24.38	& 2,079,521.71	& TL\\
\bottomrule
\end{tabular}
\caption{Results for the ESCO  Model with $|P|=2|C|$ 
%and time limit of  \red{two hours}
%one hour
\label{tab:G1}}
\end{table}

\begin{table}[tbh!]
\centering
%\scriptsize
\tabcolsep 4pt
\renewcommand
\arraystretch{1.2}
\begin{tabular}{lrrrrrrr}
\toprule			
$|C|$ & $|P|$ &  Alg  %&  \#Cuts  
&  \#Nodes 
%& $LB_{root}$ 
& DB  & Gap$[\%]$ & $\ObjVal$ &  CPU$[s]$ \\
	\cline{1-2} 
	\cline{4-8}
	%\cline{6-9}
30	&90	&com	&0	&274,839.75	&0	&274,839.75	&104.22	\\
40	&120	&com	&0	&374,154.65	&0	&374,154.65	&185.33	\\
50	&150	&com	&0	&454,088.18	&0	&454,088.18	&428.06	\\
75	&225	&com	&0	&678,652.85	&0	&678,646.43	&895.86	\\
100	&300	&com	&0	&905,017.84	&0.01	&904,970.54	&5551.27\\
150	&450	&com	&-	& -& -	&-	&TL	\\
200	&600	&com	&-	&-&-	&-	&TL	\\
\hline
30	&90	&Ben	&10566	&274,859.32	&0.01	&274,839.75	&248.72	\\
40	&120	&Ben	&98774	&387,704.33	&3.62	&374,154.65	&TL	\\
50	&150	&Ben	&41715	&481,403.55	&6.02	&454,088.18	&TL	\\
75	&225	&Ben	&65723	&722,437.01	&7.59	&671,453.24	&TL	\\
100	&300	&Ben	&30909	&947,847.23	&13.5	&835,141.79	&TL	\\
150	&450	&Ben	&31229	&1,375,809.16	&7.04	&1,285,314.13	&TL	\\
200	&600	&Ben	&615	&1,913,149.96	&13.88	&1,680,002.31	&TL	\\
\bottomrule
\end{tabular}
\caption{Results for the ESCO  Model with $|P|=3|C|$ 
%and time limit of \red{two hours}
%one hour
\label{tab:G2}}
\end{table}

\section{Conclusions}\label{sec:conclusion}
The CESs serve as an energy warehouse within each community to improve community self-consumption, therefore, clustering prosumers-based communities in smart cities and deploying the CESs with the optimal configuration is significant for more efficient energy allocation and distribution and will be an important step to realize the future decentralized electricity network and market. A paradigm shift in the community electricity network design will inevitably occur and the operational model of the power systems will also change significantly. 

In this study, we analyzed and compared the economic feasibility of CES systems in smart cities. Three %clustering 
\red{business models} are proposed to meet different stakeholders' objectives. {These objectives consider minimization of deployment/operational costs (for the island or interconnected model) or maximization of internal energy transaction profits (for the ESCO model).} Demand and generation profiles of 1,396 households from Cambridge, MA, are used as input to test the computational performance of the proposed compact models and their Benders reformulations. Moreover,  a thorough economical analysis of a neighborhood consisting of 40 consumers and 80 prosumers is conducted. For this case study, daily optimal CES operation plans are obtained and used to analyze CES performance in three different scenarios. 

Based on our simulation results, we showed that energy sharing communities deployed with CESs are more economically efficient than single households equipped with HES systems. %These benefits are even more pronounced with the increasing PV generation at prosumer locations. 
Indeed, the cost of 1.214 M\$ for deploying 80 HES at prosumer locations 
could be reduced to {551 k\$}, respectively 99 k\$, if 
%401 k\$ 
the interconnected %cluster 
model or the ESCO %cluster 
model, are employed instead. At the same time, with both models, significant savings could be achieved when it comes to the community's total electricity expenses.  Our simulations also show that, because of the high fuel consumption cost and CES with very large capacities required in the island model, this %clustering 
\red{business model} cannot be considered economical from the current technological perspective. 

{In this article we have studied deterministic mathematical models, however their accuracy can be improved by incorporating data uncertainty related to the 
%consumers' demand, prosumers' production, 
households’ load and prosumers’ generation profiles, or electricity prices. This can be achieved by, e.g., developing a two-stage stochastic programming model with a discrete set of scenarios that are sampled using sample average approximation techniques. Our empirical results demonstrate that solving the deterministic model to optimality is already a quite challenging task -- high quality solutions and small gaps can be obtained for communities with 100 to 150 households. Thus, in order to obtain computationally tractable stochastic models that can deal with data uncertainty, further sophisticated decomposition methods need to be developed. These methods can be built as an extension of the proposed Benders decomposition approach, or they can rely on alternative exact solution methods, like branch \& price, or Lagrangian decomposition. An additional extension to explore is to incorporate AC power flow constraints in our model. This results in non-linear constraints that require a customized solution method to be able to handle the non-convex problem. \red{Also as we consider operational and strategic planning, incorporating the battery's life cycle in the models is an interesting future work.} We believe that these are interesting and promising directions for future research on this topic.  }

\red{In this work, we consider that the community is sharing infrastructure and resources in terms of shared ownership of the energy storage. However, the potential collaboration to manage and optimize the energy resources collectively such as demand-side management and demand-response programs is not incorporated in the mathematical models. Thus the collective decision-making in terms of social and collaborative aspects of resource sharing is worthy of investigation as part of future work. A future direction is also to include fairness especially if some households discharge electricity more than others. As a result, depending on the collaborative nature of the community and the power consumption profiles of the households within that community, one might incorporate additional costs that are paid by the consumers for discharging from the battery. Finally, in our current setup, the number and the locations of prosumers are fixed and are part of the input. We can extend this work to also decide on the number of prosumers by running our models with different prosumers scenarios. This type of simulation can provide insights into the number and the selection of the prosumers and their impact on the system.}

%Future research should also include the combination of different business models in ESCO model. Besides, other more flexible pricing schemes should also be considered for prosumers. Finally, further research should be conducted in the optimization of CES shares allocation among prosumers.

%\newpage
\section*{Acknowledgements}
 {This research is partially supported by Natural Sciences and Engineering Research Council of Canada Discovery Grant 2017-04185 and by France Canada Research Fund. Bissan Ghaddar's research is additionally supported by the David G. Burgoyne Faculty Fellowship. The authors would like to thank the anonymous referees for their helpful and insightful suggestions for improving the paper.}
\bibliographystyle{apa}
\bibliography{ref}
%\section*{Appendix}

%\section*{Nomenclature}

%\subsection*{Abbreviations (in alphabetic order)}
%\begin{table}[H]
%\begin{tabular}{ll}
%CES  & Community Energy Storage         \\
%CN& Interconnected Clustering Model     \\ 
%DERs & Distributed Energy Resources     \\
%DSM  & Demand Side Management           \\
%EMS  & Energy Management System         \\
%ESCO & Energy Service Companies         \\
%ESS & Energy Storage System\\
%HES  & Household Energy Storage         \\
%IL & Island Clustering Model\\
%MILP & Mixed-Integer Linear Programming \\
%NPV & Net Present Value     \\
%PEV  & Plug-in Electric Vehicle         \\
%PV& Photovoltaic\\
%SOC & State of Charge
%\end{tabular}
%\end{table}

\end{document}